\documentclass[12pt,twoside]{article}
\usepackage[arrow,matrix,curve]{xy}
\CompileMatrices
\usepackage{super2}
\let\Bbb\mathbb

\expandafter\ifx\csname xymatrix\endcsname\relax%
\warndef\ar#1 #2,#3 #4,#5;{\put(#1,#2){\vector(#3,#4){#5}}}
\warndef\bx#1 #2,#3;{\put(#1,#2){\makebox[0pt]{$#3$}}}
\warndef\ln#1 #2,#3 #4,#5;{\put(#1,#2){\line(#3,#4){#5}}}
\warndef\mar#1 #2,#3 #4,#5;#6 #7,#8;{\multiput(#1,#2)%
(#6,#7){#8}{\vector(#3,#4){#5}}}
\warndef\mln#1 #2,#3 #4,#5;#6 #7,#8;{\multiput(#1,#2)%
(#6,#7){#8}{\line(#3,#4){#5}}}
\fi

\warndef\eqnum#1{\eqno (#1)}
\warndef\fnote#1{\footnote}
\warndef\eopp{\hbox{\vrule width 8pt height 8pt depth 0pt}}
\warndef\eop{\hbox{\vrule width 6pt height 6pt depth 0pt}}
\warndef\arr{\longrightarrow}
\warndef\harr{\hookrightarrow}
\warndef\ol{\overline}
\warndef\ul{\underline}
\warndef\bw{\bigwedge}
\warndef\bv{\bigvee}
\warndef\la{\langle}
\warndef\ra{\rangle}
\warndef\grot{Grothendieck\mxspace}
\warndef\grots{Grothendieck schemes\mxspace}
\warndef\tdt{\times\dots\times}
\warndef\cdc{,\dots,}
\warndef\ttilde{\widetilde}
\warndef\hhat{\widehat}
\warndef\comp{\mbox{\scriptsize$\circ$}}
\warndef\mult{\centerdot}

\rmname{mod}

\mathdef\One$\mbox{{\rm1\kern-0.37em1}}$
\mathdef\p$\mathfrak{p}$
\mathdef\Ga$\Gamma$
\mathdef\La$\Lambda$
\mathdef\Lai$\La_i$
\mathdef\Li#1$\La_{#1}$
\mathdef\Si$\Sigma$
\mathdef\set$\bf Set$
\mathdef\SET$\bf SET$
\mathdef\top$\bf Top$
\mathdef\man$\bf Man$
\mathdef\bman$\bf BMan$
\mathdef\sman$\bf SMan$
\mathdef\sreg$\bf SReg$
\mathdef\vbun$\bf VBun$
\mathdef\setg$\set^\gr$
\mathdef\mang$\man^\gr$
\mathdef\topg$\top^\gr$
\mathdef\oR$\ol\Ree$
\mathdef\oC$\ol\Cee$
\mathdef\oK$\ol\Kee$
\mathdef\oV$\ol V$
\mathdef\oVi$\ol{V'}$
\mathdef\oVii$\ol{V''}$
\mathdef\oViii$\ol{V'''}$
\mathdef\Sn$\mathfrak S_n$
\mathdef\gr$\bf Gr$
\mathdef\glutu${\bf Glut}_\cU$
\mathdef\psitu${\bf Psite}_\cU$
\mathdef\otau$\cO_\tau$
\mathdef\ox$\cO/X$
\mathdef\tilc$\widetilde\cC$
\mathdef\tild$\widetilde\cD$
\mathdef\tauo$\tau_\cO$
\mathdef\id#1${\rm Id}_{#1}$
\mathdef\K$\Bbbk$
\mathdef\Zplus$\Bbb Z_+$
\mathdef\N$\Bbb N$
\mathdef\R$\Bbb R$
\mathdef\C$\Bbb C$
\mathdef\Z$\Bbb Z$
\mathdef\Zz$\Bbb Z_2$
\mathdef\Zk$\Bbb Z_2^k$
\mathdef\nN$n\in\N$
\mathdef\iN$i\in\N$
\mathdef\iZ$i\in\Z$
\mathdef\nZ$n\in\Z$
\mathdef\0$\bar0$
\mathdef\1$\bar1$
\mathdef\implies$\ \Longrightarrow\ $
\mathdef\implieby$\ \Longlefttarrow\ $
\mathdef\lequiv$\ \Longleftrightarrow\ $

\mathdef\darr#1#2#3#4$#1,#2\colon#3{}_{\arr}^{\arr}#4$
\mathdef\sar#1,#2,#3;$#1\stackrel{#2}{\arr}#3$
\mathdef\spar#1.#2.#3;$#1\stackrel{#2}{\arr}#3$
\mathdef\fam#1#2$\{#1\}_{#2}$
\mathdef\fami#1$\fam{#1}{i\in\N}$
\mathdef\famn#1$\fam{#1}{n\in\N}$
\mathdef\faml#1$\fam{#1}{\La\in\gr}$
\mathdef\Mod#1#2${\bf Mod}_#1(#2)$
\mathdef\SMod#1#2${\bf SMod}_#1(#2)$
\mathdef\L#1#2#3#4#5${{\rm L}_#1^#2(#3_1\cdc #4_#2;#5)}$
\mathdef\LL#1#2#3#4#5$\cL_#1^#2(#3_1\cdc #4_#2;#5)$
\mathdef\T#1#2#3#4$#3_1\otimes_#1\dots\otimes_#1#4_#2$
\mathdef\diff#1${\rm Diff}(#1)$
\mathdef\sdiff#1$\hhat{{\rm Sdiff}}(\mathcal{#1})$
\mathdef\hom#1#2#3${\rm Hom}_{#1}(#2,#3)$
\mathdef\homk#1#2$\hom\K{#1}{#2}$
\mathdef\homc#1#2$\hom\C{#1}{#2}$
\mathdef\homr#1#2$\hom\R{#1}{#2}$

\usepackage{smart}
\usepackage{amsbsy}
\Aiv
\plainpage
\DropZerosNumberingStyle
\theorems[subsection]
\def\wt{\widetilde}
\def\ol{\overline}
\def\ul{\underline}
\def\iI{_{i\in I}}
\def\uU{_{u\in U}}
\def\la{\langle}
\def\ra{\rangle}
\rmname{Sh}
\oplim{colim}
\mathdef\Tsub$T_{sub}$
\mathdef\Csub$\cC_{sub}$
\mathdef\Osub$\cO_{sub}$
\mathdef\tsub$\tau_{sub}$
\mathdef\vd$\vdash$
\mathdef\vD$\vDash$
\mathdef\Shu$\Sh_{\cU'}\cC$

\def\version{1}
\def\vers#1#2#3{#1}
\begin{document}
\title{Glutoses: a Generalization of Topos theory}
\author{Vladimir MOLOTKOV\thanks{
\noindent{\sl 1991 Mathematics Subject Classification} 18F10 (Primary) 18F15,18F20,14A20,14A22 (Secondary)\hfill\hfill\protect\linebreak
\indent{\sl Key words and phrases:} Category
theory, topos theory, algebraic geometry.\hfill\hfill\protect\linebreak
\indent Research partially
supported by the Ministry of Sci. and Educ. of Bulgaria
grant F-610/98-99.}\\
Inst. for Nuclear Res. and Nucl. Energy\\
blvd.\ Trakia 72, Sofia 1784, Bulgaria\\
e-mail: {\sl vmolot@inrne.bas.bg}}
\date{December 8, 1999}%
\maketitle

\begin{abstract}
A generalization of topos theory is proposed
giving an abstract realization of such categories
as, say, the categories of manifolds and of
Grothendieck schemes on the one hand, and permitting one,
on the other hand, a view on ``non-commutative'' or,
more generally, ``universal'' algebraic geometry,
which is alternative to already existing, and is closer,
in some sense, to the classical
Grothendieck's construction of commutative schemes.
Another immediate application of the theory developed
is construction of an extension of the category of Grothendieck
schemes to the category of ``\'etale schemes''
containing together with any scheme every \'etale sheaf
over it as well.

The main result of this work is that for any presite
satisfying some smallness conditions (existence of {\it local\/} sets of
topological generators) there exists the universal ``completion'' of a presite
to a glutos.

This paper is a corrected and extended version of
JINR-preprint E5-93-45 (1993) \cite{Mol1}.
It is more or less complete
as concerns formulations of results, but proofs are only hinted in several
places. A version which
includes proofs as well is now in progress.
\end{abstract}

\section{History \& Motivations}
\label{history}
The theory presented here originates from autor's works on the theory of
infinite-dimensional supermanifolds (\cite{MolS1}-\cite{MolS3}).
It was purposed originally just to develop the technical tools permitting
one to deal automatically with numerous kinds of ``charts and atlases''
arising in the theory of supermanifolds (in the latter theory the role of
the category of sets is played by some functor category, equipped with
some pretopology).

The scheme is as follows: given a category \cC equipped with some
{\it pre\/}topology (which has {\it not\/}, generally speaking, a set
of topological generators) to construct another category $\wt{\cC}$
with a pretopology together with the universal continuous functor
\begin{equation}
\label{0}
Y_{\cC}\colon\cC\harr\wt{\cC}
\end{equation}
such that, roughly speaking, $\wt{\cC}$ is complete with respect
to ``gluing along open ($=$ belonging to some covering) arrows'' and,
besides, its objects are ``locally isomorphic'' to objects of \cC.

A bit more precisely, this means that certain functors with values
in open arrows of $\wt{\cC}$ (called `glueing data' or `gluons')
have colimits satisfying certain conditions. I have called such
universal completion `glutos associated with \cC'.

An archetypical example is the case when \cC is the category
of open regions of Euclidean spaces (resp. of Banach spaces, resp.
of Banach superspaces) with (super)smooth morphisms.
Then $\wt{\cC}$ is naturally equivalent to the category \man
of smooth manifolds (resp. to the category \bman of smooth Banach manifolds,
resp. to the category \sman of Banach supermanifolds), whereas any glueing
data can be identified with some atlas on its colimit.

Essential point here is that the glutos $\wt{\cC}$
is determined completely
(up to natural equivalence, of course) by the category \cC and
its pretopology.

Originally there were wery strong restrictions imposed on the pretopology
of \cC in order that the construction of $\wt{\cC}$ by means of
``charts and atlases'' routine make sense. In particular:
a) the pretopology on \cC must be subcanonical;
b) open arrows should be mono;
c) union of open arrows must exist and be open again (as is the case
in the ``generic'' pretopology on \top);
d) open subobjects of any object of \cC must form a set.

Later on the construction of charts and atlaces was extended, so as to include
pretopologies not satisfying condition c) above. This is, e.g., the case
when \cC is the category dual to that of commutative rings, with
Zariski pretopology on it. It turns out that in this case the category
$\wt{\cC}$ is naturally equivalent to the category of \grots.

The latter result forces a natural question: whether one can
remove condition a) on the pretopology of \cC to construct
{\it non-commutative\/} schemes by means of the ``glutos generator''~(\ref{0})
(there are several analogues of Zariski pretopology for the category
dual to the category of {\it all\/} rings (see~\cite{Coh}), neither of them
is subcanonical), as well as to get rid of condition b), to be able
to produce the category of algebraic spaces (or, may be some its extension)
by means of the same universal glutos construction.

The answer to the first of this questions is positive.

It turned out that, in case of a pretopology satisfying conditions a)
and d) above, the universal arrow~(\ref{0}) exists and can be represented
as the composition of three universal arrows:
\begin{equation}
\label{0'}
Y_{\cC}\colon
\sar\cC,Y'_{\cC},{\cC_{sub}};
\sar{},Y''_{\cC},\ol{\cC};
\sar{},Y'''_{\cC},\wt{\cC};,
\end{equation}
where the functor $Y'_{\cC}$ is the universal arrow 
into presites with subcanonical pretopology (cf.~Sect.~\ref{Proof}),
the functor $Y''_{\cC}$ is the universal completion of $\cC_{sub}$
by objects which are ``unions of open subobjects'', whereas the
functor $Y''_{\cC}$ is the ``charts and atlases'' construction in
``canonical'' sense.

The scetch of this early results was published in 1986 in~\cite{Mol}.
The manuscript containing detailed results and its proofs
(named ``Manifolds IV'') was never published.%
\footnote
{I am very grateful to D.Leites, who organized translation of it
into ChiWriter in 1988
(in those times \TeX\ was not yet, unfortunately, as popular as it is now.).
I hope the corresponding ChiWriter file,
translated by Chi2TeX programm (though badly: Chi2TeX knows nothing about
commutative diagrams!),
will reduce essentially the work needed to write the complete
version of this work on glutoses.}

During next several years I have made a number of attempts
to generalize the main result of~\cite{Mol}, trying, in particular,
to get rid of the condition b) (open arrows are mono)
on the pretopology of \cC while {\it extending\/}, simultaneously,
the very notion of a
glutos. All this attempts were non-satisfactory,
because, as I understand now, I was too hypnotized by ``charts
and atlases'' paradigm. It turned out now that, in the general case,
there are much more new objects in the glutos $\wt{\cC}$, then
can be obtained by using atlases alone (i.e. glueing along open
arrows in \cC only). The whole process of completion of \cC
to $\wt{\cC}$ turns out to be transfinite!

Meanwhile, in 1993 there appeared a paper of Paul Fait~\cite{F} devoted
to the same problem, whose main result (Theorem 14.4 on p.94)
was the proof of existence of the completion
$Y_{\cC}\colon\cC\harr\wt{\cC}$
for a presite \cC with subcanonical pretopology but
at formally weaker other conditions imposed on the pretopology of \cC
than those in ~\cite{Mol} (they reduce exactly to conditions
in~\cite{Mol}, plus the requirement that the pretopology of \cC be
subcanonical, if one assumes that open arrows are mono).

Unfortunately, the work~\cite{F} uses highly non-standard
terms even for quite standard categorical notions,
which makes it difficult for reading. Besides,
just to understand what the main theorem 14.4 of~\cite{F}
is talking about, one needs to look through the whole
paper~\cite{F} in search of definitions.

That's why 
below is given the translation (to a reasonable extent) of the main result
of~\cite{F} to common math language purposing to compare
this result with
that
of~\cite{Mol}.

A category \cC equipped with a subcanonical (``intrinsic''
in terms of P.~Fait) pretopology $\tau$ is called by P.~Fait
a {\bf local structure} if
\begin{list}{}{\topsep1mm\itemsep1mm}
\item[a)] for any object $X$ of \cC equivalence classes of open arrows
$u\colon U\arr X$ form a small set (smallness condition);
\item[b)] Any family of open arrows such that there exists
a refinement of it which is a covering of $\tau$ is a covering itself
(such pretopology is called ``flush'' by P.Fait; in the present paper
the corresponding property is just added to the definition
of a pretopology --- see condition (PT4) in Sec.~\ref{ps});
\item[c)] Every open retraction $u$ is a covering morphism
(i.e. the singleton $\{u\}$ is a covering);
\item[d1)] Every arrow $u\colon U\arr X$ which is a covering locally
is a covering of $X$. Here ``$u$ is a covering locally''
means that there exists a covering
$\{u_i\colon U_i\arr X\}_{i\in I}$ of $X$
such that for any $i\in I$
the pullback projection $U_i\prod_XU\arr U_i$ is a covering of $U_i$;
\item[d2)] Let $u\colon U\arr X$ be an arrow such that
there exists a covering
$\{u_i\colon U_i\arr U\}_{i\in I}$
such that for any $i\in I$
the restriction arrow $uu_i\colon U_i\arr X$ is open and, besides,
the pullback projection $U_i\prod_XU\arr U_i$ is a covering of $U_i$.
Then $u$ is an open arrow.

(conditions d1)-d2) together is the ``CLCS condition'' in the terminology of P.Fait);

\item[e)] For the last condition of P.~Fait: ``\cC is complete under Cvm''
I failed to find a simple formulation in terms of ``internal'' properties
of \cC itself; as formulated by P.~Fait it includes a complicated construct,
namely, some functor
$+\colon\cC\arr\cC^+$
called by P.~Fait ``the plus functor''. The construction
of this ``plus functor'' takes a considerable part of the paper~\cite{F}.
First, there is constructed the presite $\cC^p$, whose objects are
just glueing data on the presite \cC (``canopies'' in the terminology
of P.~Fait) and whose arrows and pretopology are ``forced'' in some
sense by that of \cC. The pretopology on $\cC^p$ fails to be subcanonical,
so the next step is to form a new
category $\cC^+$ out of $\cC^p$, this time with subcanonical pretopology by
some process, called ``smoothing'' by P.~Fait. I suspect (though have not
checked it) that ``smoothing'' is equivalent to a particular case of
our construction of universal presite with subcanonical pretopology
(see Section~\ref{Proof} below).
Very roughly and informally the ``plus functor'' turns out to be the
``square root'' of the desired functor
$Y_{\cC}\colon\cC\harr\wt{\cC}$.

In terms of this ``plus functor'' the last condition in the definition
of local structure says that the functor $+$ reflects covering morphisms,
in the sense that a single arrow $u\colon U\arr +(X)$ of $\cC^+$
is a covering of $+(X)$ iff it is isomorphic (in $\cC^+/+(X)$) to an
arrow $+(u')$ for some arrow $u'\colon U'\arr X$ of \cC.
\end{list}

A local structure \cC is called by P.~Fait a {\bf global structure} if
it ``is complete under affinization'', which means in our terms
that every open glueing data in \cC (see definition in Sec.~\ref{glu} below)
has a universal colimit satisfying certain natural conditions
(strong form of our conditions {\bf(${\bf G4_{\cU}+5}$)}
in Sect.~\ref{glu} below).

The p.~(E) of the main theorem~14.4 of~\cite{F} states then that
{\it for any local structure \cC the category $\cC^{++}$ is a global
structure, whereas the functor $++\colon\cC\arr\cC^{++}$ is universal
among continuous functors to global structures\/}.

One can see that if one imposes on a local structure
the additional requirement that all open arrows are mono, then
conditions c)--e) in the definition of the local structure
above can be omitted. In this particular case `local structure'
resp. `global structure' of P.~Fait is exactly the same thing as
`preglutos with subcanonical topology' resp. `glutos' in my paper~\cite{Mol}
(in the current paper glutoses in the sense of~\cite{Mol} are called
`nearly \cU-glutoses', and they correspond to a particular
case of glutoses in the present sense, namely, SG-glutoses).

On the other hand, without the requirement that opens are mono,
the condition e) in the definition of local structure
seems to me to be too {\it ad hoc} and, besides, I see
no simple way to check whether a given presite satisfies
it. For example, the category of algebraic spaces with
\'etale pretopology on it is declared to be
a global structure (Example~E on p.2 of~\cite{F}),
but I have found no indication in the text of~\cite{F}
why the condition e) is valid for this case.

In short, it is not clear at all, to what extent
conditions imposed on the pretopology of a local
structure are weaker compared to the same conditions
plus the condition stating that opens are mono.

In any case, the proof in~\cite{F} of existence of the universal
global structure is performed completely inside a given universe,
without using the axiom of existence of strongly inaccessible
cardinals, whereas I use systematically in proofs ``big'' toposes
(not belonging, generally speaking, to the universe where `local'
and `global' structures live), which makes the results of~\cite{F}
stronger than mine from purely set-theoretic point of view.

My own attempts to generalize the results of my paper~\cite{Mol}
were based on the following
experimental fact (discovered after the publication of~\cite{Mol}):
for any object $X$
of a `glutos' \cC the full subcategory of the category $\cC/X$
consisting of {\it locally open\/} arrows $f\colon Y\arr X$
turnes out to be a {\it topos}. This generalizes the (archetypical!) fact
that the topos of sheaves over a topological space $X$
is naturally equivalent to the category of local homeomorphisms over $X$.

Toposes arising in such a way are of a very special kind --- they are so
called SG-toposes (i.e. those toposes, for which subobjects of $1$ form
a set of generators). So it seems to be natural to look
for such an extension of the theory, where `locally open' arrows over
any object of a `generalized glutos' form a topos again, this time
{\it arbitrary} \grot topos.

It took me several years both to find the right generalization
of glutos definition and to get rid of {\it all} restrictions
on the pretopology of \cC in the ``universal glutos generator''~(\ref{0}),
excepting a highly weakened form of smalness condition d) above
(the existence, for any object $X$ of \cC,
of a small set of topological generators in the induced presite of
open arrows over $X$).

The axioms of the theory arising admit an elementary version,
in which case the theory may be viewed as some generalization
of the theory of elementary toposes (see Sect.~\ref{def}), whereas the version
of its axioms including small (with respect to some universe \cU)
families of objects or arrows in its formulation
(see Sect.~\ref{ug}) is a generalization of \grot's toposes.

And it turned out that in the general case charts and atlases are
non-adequate for constructing the glutos $\wt{\cC}$ out of \cC.

Unfortunately, the general definition of a glutos given in
the first version of this paper~\cite{Mol1}
turned out to be wrong. It failed to satisfy to what physicists
would call ``the correspondence principle'':
the old theory must be a particular case of the new one,
meaning here that the 2-category of ``old'' glutoses
(i.e. in the sense of~\cite{Mol})
must be (naturally equivalent to) a full 2-subcategory of ``new'' ones
(i.e. in the sense of~\cite{Mol1}).

Namely,
one of the axioms of glutoses
(a part of the axiom (G5)) as formulated
in~\cite{Mol1}
turned out to be too strong, so that all classical examples
(topological spaces, manifolds, etc.) failed to satisfy it%
\footnote{The counterexample was found by
Prof.~P.~T.~Johnstone (1983, private communication);
I am very grateful to him
for discovering this nasty bug in my definitions,
leading to bugs
in a number of statements
in the paper~\cite{Mol1})
},
though ``glutos generator'' above works for this definition, too!

It took me only about a month to correct the definition,
finding a weaker theory, satisfying the correspondence principle,
but the writing of this corrected version of~\cite{Mol1}
was delayed for more than 5-years due to various reasons of
non-mathematical nature.

As to the ``wrong'' glutoses, they are not thrown away.
Instead they are present in this paper under the name
`ultraglutoses'.
Though they do not satisfy correspondence principle
for old version of glutos theory, scetched in~\cite{Mol}, they
{\it do satisfy} this principle, being one of the correct extensions
of the {\it topos theory\/} instead!

Glutoses do satisfy the latter principle as well,
but, in some sense, the theory of ultraglutoses is a better
extension of the topos theory, because they are more ``ideal''
from purely categorical point of view (more diagrams have
colimits inside ultraglutoses and this colimits are ``good'').
The most important thing is that any glutos has the universal
fully faithful imbedding into an ultraglutos.

But ultraglutoses themselves seem to be not the end of the story yet.

For example, the axiom of finite completeness is not included
in their definition, just because I wrongly believed originally
that such
categories as \man, etc.
are to be models of this theory.

Now, there are no principal
obstructions to add the finite completeness axiom to the theory
of ultraglutoses, as far as one proves that {\it any glutos \cC imbeds
{\rm(fully and faithfully)}
into finitely complete ultraglutos $\cC'$, so that the imbedding
is universal among all continuous functors {\rm(not necessarily exact!)}
into ultraglutos\/}.
I believe the latter statement is true,
but have not find time yet to check it.

\section{Set-theoretic conventions}
\label{sets}
We will work here within Morse set theory (see \cite{Mor}),
with the usual axiom added stating the existence for any
set of a universal {\it set\/} containing it.

The latter axiom seems to be not necessary for the
validity of our main Theorem~\ref{th} (reformulated properly),
but its use highly facilitates our construction of universal glutoses
and ultraglutoses.

Namely, to constuct a universal (ultra)glutos $\wt\cC$ for a presite \cC
living in some universe \cU, we use as building blocks
objects and arrows of the topos \Shu of shieves on \cC
with values in some {\it higher} universe $\cU'$, containing \cU
as an element. This simplifies a number of proofs.

All categories, presites, toposes, glutoses, etc. are
supposed to be sets, so that there exists the legitimate
2-category of all categories, resp. presites, etc.
(which is a proper class); similarly, (pseudo)functorial
operations defined in several places below on objects,
arrows and 2-arrows of 2-categories above are incorporating
together to produce legitimate 2-(pseudo)functors.

Here are assumed  the definitions of \cite{Mor}
for ordered pairs and, more generally, families (= tuples
in the terminology of \cite{Mor}) so that for any universe
\cU (including the biggest universe of {\it all\/} sets)
a family of subclasses of the universe \cU indexed by a
subclass of the universe \cU is again a subclass of \cU
and behaves well%
\footnote{This means that families are ``separated'':
$\{\cC_i\}_{i\in\cI}=\{\cC'_i\}_{i\in\cI'}$
iff $\cI=\cI'$ and for any $i\in\cI$ one has $\cC_i=\cC'_i$.
In other words, there is no loss of information while encoding
some data inside a family, even the big one.}.
In fact, it is just this choice of
definition for families which permits one to {\it define\/}
the ``big'' 2-categories above and 2-(pseudo)functors
between them as terms of Morse set theory: e.g.,
a 2-category \cC {\it is\/} a finite tuple
$\langle\cC_0,\cC_1,\cC_2,\dots\rangle$ of classes
satisfying certain conditions.

\section{Inverse images of equivalence relations}
\label{eqrel}
Here are given some definitions and elementary results, necessary
to formulate the axioms of glutoses.

\begin{Prop}
Let
$
d_0,d_1\colon U{}_{\arr}^{\arr}X
$
be an equivalence relation in a category \cC (not necessarily
with finite products). Let $f\colon X'\arr X$ be a pullbackable arrow
of \cC. Consider the diagram of three pullbacks

\begin{equation}
\label{indeqrel}
\expandafter\ifx\csname xymatrix\endcsname\relax%
\divide\dgARROWLENGTH by2
\begin{diagram}
\node{}\node{}\node{U'}\arrow{se}\arrow{sw}\\
\node{}\node{U'_0}\arrow{se}\arrow{sw}\node{}\node{U'_1}\arrow{se}\arrow{sw}\\
\node{X'}\arrow{se,t}{f}\node{}\node{U}\arrow{se,t}{d_1}\arrow{sw,t}{d_0}%
\node{}\node{X'}\arrow{sw,t}{f}\\
\node{}\node{X}\node{}\node{X}
\end{diagram}
\else
\vcenter{%
\xymatrix@=10pt{
&&U'\ar[dr]\ar[dl]\\
&U'_0\ar[dr]\ar[dl]&&{U'_1}\ar[dr]\ar[dl]\\
{X'}\ar[dr]^f&&U\ar[dr]^{d_1}\ar[dl]_{d_0}&&X'\ar[dl]_f\\
&X&&X
}}
\fi
\end{equation}
Let $d'_0\colon U'\arr X'$ {\rm (resp. $d'_1\colon U'\arr X'$)} be
the composition of
left
{\rm (resp. {\it
right\/
})} pullback
projections of the diagram above.
Then the pair
$
d'_0,d'_1\colon U'{}_{\arr}^{\arr}X'
$
is an equivalence relation on $X'$ which will be called {\bf induced on
$X'$ by $(d_0,d_1)$ along $f$} or {\bf inverse image of $(d_0,d_1)$
along $f$} and will be denoted $f^*(d_0,d_1)$ {\rm (Note that
$f^*(d_0,d_1)\not=(f^*d_0,f^*d_1)$!)}.
\end{Prop}

\noindent{\bf Proof.} Imbed \cC into some ``big'' topos
$
\widehat{\cC}=[\cC^\circ,{\bf Set}_{\cU}]
$
of preshieves on \cC with values in a universe \cU via Yoneda
(recall that due to our choice of set theory such universe exists if
the category \cC is a set). Augmenting the diagram (\ref{indeqrel})
by a coequalizer $r\colon X\arr R$ of $(d_0,d_1)$ in $\widehat{\cC}$
and using Giraud theorem as well as left exactness of Yoneda functor
we conclude that $(d'_0,d'_1)$ is a kernel pair (in $\widehat{\cC}$)
of $rf$ which implies that $(d'_0,d'_1)$ is an equivalence relation.\ \eop

\section{Glutoses: definition}
\label{def}
An {\bf elementary  glutos} is a kind of ``generalized elementary topos''
(not to confuse with quasitoposes!): it is a category \cC, equipped with a
suitable structure, given by a subset \cO of arrows of \cC  (elements of
which will be called {\bf open\/} arrows of \cC), which generates over any
object $X$ of \cC an elementary topos, namely, $\cO/X$ . Here $\cO/X$
denotes the full subcategory of $\cC/X$ formed by all objects which are
arrows of \cO. So that, in some sense, a glutos is a topos {\em locally\/}
(not to be confused with local toposes!). If one grasps, metaphorically,
a glutos  as a family of toposes coherently glued together into a single
category \cC by means of a ``glueing'' structure \cO,
then the term `glutos' itself can be thought of as an
abbreviation for `GLUed bunch of TOposeS'. An alternative interpretation of
this term: in glutoses one can {\it glue\/}
together finite families of objects along open arrows (see section \ref{glu}
below for details).

Exact conditions \cO must satisfy (``axioms of elementary glutoses''), are
the following conditions {\bf(G1)}-{\bf(G5P)} below.

But we need first to give some definitions (to formulate correctly
condition {\bf(G5c)}). An equivalence relation
$u,v\colon U{}_{\arr}^{\arr}X$ will be said to
be \cO-{\bf coequalizable}, if there exists an arrow $q\colon X\arr Q$
belonging to \cO such that $qu=qv$; the relation
$u,v\colon U{}_{\arr}^{\arr}X$ will be said to
be {\bf (finitely) locally \cO-coequalizable} if there exists a (finite)
epi family $\{u_i\colon U_i\arr X\}\iI$ of arrows of \cO such that
every induced equivalence relation $u_i^*(u,v)$ is \cO coequalizable
(in the latter definition it is supposed that both $u$ and $v$ are
pullbackable so that the corresponding induced equivalence relations do exist).

\begin{list}{}{\topsep1mm\itemsep1mm}
\item[\bf(G1)] \cO contains all iso's of \cC, is contained in the
set of all pullbackable arrows of \cC and is stable by composition and
pullbacks;

\item[\bf(G2)] $fg\in \cO$ and $f\in \cO$ implies $g\in \cO$;

\item[\bf(G3)] a) For any object $X$ of \cC the category $\cO/X$ is
an elementary topos and

b) for any $f\colon X\arr Y$ in \cC the
functor $f^{*}\colon\cO/Y \arr\cO/X$ (which is defined,
due to (G1), and is left exact) is an inverse image of some geometric
morphism;

\item[\bf(G4)] a) \cC has disjoint and universal finite coproducts, such
that canonical injection morphisms belong to \cO;

b) for any finite family
$\{U_i\stackrel{u_i}{\arr}X\}\iI$ of arrows of \cO the
colimit arrow
$\coprod_i U_i \arr X$ belongs to \cO.



\item[\bf(G5)] a) any epi of \cC, which belongs to \cO, is effective;

\noindent b) if both $fp$ and $p$ belong to \cO and $p$ is epi
then $f$ belongs to \cO;

\noindent c) Any equivalence relation
$u,v\colon U{}_{\arr}^{\arr}X$ in \cC
which is {\bf open} {\rm(i.e. $u,v\in \cO$)} and finitely locally
\cO-coequalizable, is effective and has a universal coequalizer in
\cC which belongs to \cO.

\item[\bf(G5P)] For any pullback diagram
\begin{equation}
\label{locpb}
\expandafter\ifx\csname xymatrix\endcsname\relax%
\divide\dgARROWLENGTH by2
\begin{diagram}
\node{V}\arrow[2]{e}\arrow{s}\node[2]{Y}\arrow{s,r}{f}\\
\node{X}\arrow{e,t,A}{r}\node{R}\arrow{e,t}{g}\node{S}
\end{diagram}
\else
\xymatrix{
V\ar[rr]\ar[d]&&Y\ar[d]^f\\
X\ar@{>>}[r]^r&R\ar[r]^g&S
}
\fi
\end{equation}
such that $r$ is an open epi
there exists a pullback of $f$ along $g$.

\end{list}

A glutos $(\cC,\cO)$ will be called an {\bf ultraglutos} if
p.~c) of the axiom (G5) is replaced by the stronger axiom:

\begin{list}{}{\topsep1mm\itemsep1mm}

\item[cu)] Any open equivalence relation
$u,v\colon U{}_{\arr}^{\arr}X$ in \cC
is effective and has a universal coequalizer in
\cC which belongs to \cO.

\end{list}

In what follows we will say, if necessary, `weak axiom (G5)'
or `strong axiom (G5)' to distinguish between glutoses
and ultraglutoses.

It is clear, that elementary (ultra)glutoses are really models
of some first-order theory, which is an extension of the elementary
theory of categories by some unary relation symbol \cO
($\cO(f)$ meaning `$f$ is an open arrow'), with corresponding
translations of (G1)--(G5P) added as axioms.

\begin{Rem}
\label{clop}
Structures \cO on a category \cC satisfying condition (G1) occur
so frequently that deserve, in author's opinion, to be christened
somehow. Here it is proposed to call them {\bf cloposes},
whereas for its elements is reserved the name {\bf clopen}
arrows. An argument in favour of this strange choice
of names is that in the category of topological
spaces both the class of all arrows  isomorphic to
inclusions of closed subspaces and that isomorphic
to inclusions of open subspaces satisfy condition (G1).
The terms `closed' and `open' can then be reserved
to denote something more special than simply elements of a class of
arrows satisfying condition (G1) (e.g. closed
arrows of a closure operator \cite{Kos} or open
arrows in the sense of \cite{JoMo} \footnote{I
am grateful to Prof.P.T.\ Johnstone who turned my
attention to the preprint \cite{JoMo} sending me a copy of it.}
or in the (different)
sense of the present work).

In the present work the term `open' will often be used instead
of `clopen', at least in the contexts of both glutoses and presites.
\end{Rem}

\begin{Rem}
\label{G6}
One can prove that axiom (G5P) follows from another axioms in case of
ultraglutoses (see Prop.~\ref{pull} below) as well as
in case of SG-glutoses defined in Sect.~\ref{oper} below.
\end{Rem}

%
%
%
%
%
%
%
%

\begin{Rem}
\label{pizdezh}
So we have two theories now, both pretending to be something like
``topos theory for the 2-category of cloposes''. Which one of them is
``the'' theory?

Glutos theory seems to be fixed,
because it was designed to have
as its models such creatures, existing in real Universe, as \man, \top,
etc. So that I see no way to make it stronger.

Ultraglutos theory has no such restrictions. Having no other
experimentally discovered models so far, excepting toposes themselves,
it can permit additional axioms to be added, as far as
ultraglutoses, generated via ``ultraglutos generator''~(\ref{0})
from presites (or glutoses) belonging to the Universe,
belong to the Universe as well.

One evident candidate for such an extension is
the axiom of finite completeness, as I have already
noted in Sect.~\ref{history}.

As soon as it will be proved that any ultraglutos has the universal
imbedding into complete ultraglutos, the complete ultraglutoses
may pretend that they are the ``ideal'' extension of glutos theory.

One is to stress, however, that left exactness must not be
included in the definition of morphisms between complete ultraglutoses,
as far as we will have the possibility to continue such morphisms
of glutoses as
\[ {\rm spec\colon } {\bf Schem} \arr {\bf Top}\ ,\]
to morphisms of finitely complete ultraglutoses.
\end{Rem}

\section{Examples and Counterexamples}
\label{exa}
{\bf(0)} Any topos is, canonically, an ultraglutos (set \cO=\cC);
vice versa, if a pair $(\cC,\cC)$ is a glutos,
then \cC is a topos iff it has a terminal object
(the latter condition is necessary as one can see
from example at the end of section \ref{2cat}).
For another examples of glutos structures on toposes
(\'etale structures satisfying collection axiom of \cite{JoMo})
see Remark \ref{etstr} in section~\ref{ug} below.

I know of no other ``natural'' example of ultroglutoses, though,
as our main theorem~\ref{th} shows, every glutos has a universal
full imbedding into an ultraglutos.

Archetypical examples of glutoses which are not toposes are:

\noindent {\bf(1)}~Topological spaces ({\bf Top}) with open arrows
being local homeomorphisms;

\noindent {\bf(2)}~Smooth manifolds ({\bf Man}) with open arrows being
local diffeomorphisms; for a natural number $n$ the
full subcategory ${\bf Man}_n$ of {\bf Man} consisting of all manifolds
of dimension $n$ with the empty manifold
$\emptyset$ added, with open arrows as above. The glutos
${\bf Man}_0$ degenerates, evidently, to the topos {\bf Set}
of sets, whereas ${\bf Man}_n$ for $n \neq 0$ give examples
of glutoses without terminal objects;

\noindent {\bf(3)}~Locally trivial vector bundles over smooth manifolds
({\bf Vbun}) with open arrows being just those arrows,
whose image in {\bf Man} under forgetful functor is open;

\noindent {\bf(4)}~Grothendieck schemes ({\bf Schem}) or
\mbox{$C^{\infty}$-schemes} of Dubuc \cite{Dub} ($C^{\infty}$-{\bf Schem})
with, e.g., open arrows in {\bf Schem} being morphisms which locally are
inclusions of open subschemes:
i.e., $(u\colon U \arr X)\in \cO$, if there exists
a covering of $U$ by open subschemes such that the restriction of $u$ on any
element of this covering is isomorphic to inclusion of an open
subscheme of $X$.

It is, implicitly, assumed in examples (1)--(4) above that, say,
all topological
spaces of {\bf Top} are elements of some universe
\cU, which, moreover, contains, for the cases
(2)--(3) as well as $C^{\infty}$-{\bf Schem}, the universe
$\cU_f$ of finite sets as an element. So that we will
write further, if necessary, ${\bf Top}_{\cU}$, resp.
${\bf Set}_{\cU}$, etc., instead of ${\bf Top}$, resp.
${\bf Set}$, etc..

\begin{Rem}
\label{Jcounter}
(P.~T.~Johnstone's counterexample) The fact that categories \man and \top
does not satisfy the strong form uc) of axiom (G5), i.e. are not ultraglutoses
in current terminology, was told to me by P.~T.~Johnstone in 1993.

The following proposition is a generalization of this
counterexample.
\begin{Prop}
Let $G$ be a discrete Lie group acting smoothly on a manifold $M$
via $\mu\colon G\times M\arr M$. Then:
\begin{list}{}{\itemindent-\parindent\topsep1mm\itemsep1mm}
\item[\rm a)] both $\mu$ and projection map $\pi_M\colon G\times M\arr M$
are local diffeomorphisms;
\item[\rm b)] If $G$ acts freely on $M$, then the pair $(\pi_M,\mu)$
is an equivalence relation on $M$, whose coequalizer in \top is the space
of orbits of $G$ equipped with the factor topology.
\end{list}
\end{Prop}
\end{Rem}

Of course, any smooth action of any Lie group $G$ on a manifold
lifts to the smooth action of the same group $G$ equipped with
the discrete topology.

P.~T.~Johnstone's counterexample is obtained, if one takes $M=S^1$ (circle
with unit length), $G=\Zee$, and the action of $\Zee$ on $S^1$ is generated
by shift (=rotation) on irrational distance. It is clear that the coequalizer
of the corresponding equivalence relation in \man is a single point,
and neither in \man nor in \top the coequalizer is a local homeomorphism,
hence, the strong form of axiom (G5) is not satisfied neither
in \man nor in \top.

Another evident counterexample is ``irrational'' action
of \Ree on torus $S^1\times S^1$.

So, it is a curious fact, that if one completes both glutoses \top
and \man to universal ultraglutoses (via Main Theorem~\ref{th}),
then the ``good'' orbit spaces always exist inside these completions.
Besides, in case of the ultraglutos $\wt{\man}$ the orbit space will
always have the tangent bundle as well! This just follows trivially
from our Main Theorem~\ref{th}.

And one can pose a question (though rhetorical, yet), whether one can not get
an alternative (and more simple) version of ``non-commutative
differential geometry'' of A.~Connes (or, at least, part of it),
if one develops this theory inside an ultraglutos $\wt{\man}$
in place of the glutos \man itself?

\section{The 2-category of glutoses}
\label{2cat}
A {\bf morphism} of a glutos $(\cC,\cO)$ into a glutos
$(\cC',\cO')$ is a functor $F\colon\cC \arr
\cC'$ between underlying categories, which respects the structures
involved, i.e. satisfies the following conditions
below:

\noindent {\bf(MG1)}~$F(\cO)\subset \cO'$ and $F$ respects
pullbacks of open arrows along arbitrary arrows of~\cC;

\noindent {\bf(MG2)}~For any object $X$ of \cC the induced
functor $F/X\colon\cO/X \arr\cO'/FX$ (which
is defined and is left exact, due to (MG1)) is an inverse image of some
geometric morphism.

Examples of morphisms of glutoses are: the string of
forgetful functors
\[ {\bf VBun \arr Man \arr
Top \arr Set};\]
the tangent functor
$T\colon {\bf Man}\arr {\bf VBun}$, as well as
the transition to base manifold functor
$B\colon {\bf VBun}\arr {\bf Man}$; the functor
\[ {\rm spec\colon } {\bf Schem} \arr {\bf Top}.\]

Besides, all natural functors ``inside'' a glutos
are morphisms of glutoses as shows the following

\begin{Prop}
\label{canmgl}
For any object $X$ of a glutos $(\cC,\cO)$,
the source functor $d_0\colon $\cO/X$ \arr
\cC$ is a morphism of glutoses; for any arrow
$f\colon X \arr Y$ of \cC the functor
$f^*\colon\cO/Y \arr\cO/X$ is a morphism
of glutoses.
\end{Prop}

Of course, the toposes $\cO/X$ and $\cO/Y$
above are considered as glutoses via canonical
glutos structure of example (0) of Sect.\ref{exa}.

Adding any natural transformations between morphisms of
glutoses as 2-arrows
one obtains a 2-category {\bf Glut} of glutoses
as well as its full subcategory {\bf UGlut} consisting of
ultraglutoses.

Now, as is clear enough, the 2-category of toposes
imbeds contravariantly
to that of glutoses ``almost fully''
in the sense that any morphism of glutoses
$f\colon\cE\arr\cE'$ between
toposes $\cE$ and $\cE'$ decomposes as
\[ \cE\approx\cE/1\stackrel{f/1}{\arr}
\cE'/f1\stackrel{d_0}{\arr}\cE'\ ,\]
which means that the ``deviation''
of a glutos morphism~$f$
between toposes from an inverse image of some geometric
morphism is just the difference between $f1$ and $1$;
the imbedding above would be full if one permits for inverse images
of geometric morphisms not to respect terminal objects.

Nevertheless,
the theory arising is not just generalization of topos theory but, rather,
a counterpart to the latter. An essential difference is that {\bf presites\/}
(=categories
\cC, equipped with a pretopology $\tau$
\footnote{Note only that we will consider sinks in \cC and, in particular,
coverings of $\tau$ as elements of the set ${\rm Ob}\cC
\times \cP({\rm Mor}\cC)$ (where \cP stands for
power set), rather than indexed families of arrows of \cC,
though in practice indexed families will be used as well, as representing,
in an evident way, ``real'' sinks. The set of all pretopologies on \cC
forms then a closure system (in the sense of \cite{Coh1}) on the set of all
pullbackable sinks of \cC.}) play for glutoses the same role sites play
for toposes, as will be seen in section \ref{ugps} below%
\footnote{\grot topologies are, clearly, non-adequate, because there
is no natural subset \cO of ``open'' arrows associated with them.
In Appendix~\ref{ggt} is suggested a generalization of both clopos
definition and that of \grot topology ({\it depending} on generalized clopos
structure), which, I beleive, will permit to prove a stronger form of
the Main Theorem~\ref{th} by using generalized \grot topologies
in place of pretopologies.
}.


We conclude this section defining a {\bf subglutos} of
a glutos $(\cC,\cO)$ as a glutos $(\cC',\cO')$
such that $\cC'$ is a subcategory of \cC closed
with respect to composition with isomorphisms of \cC,
the set $\cO'$ is a subset of \cO and, besides, the
inclusion functor $\cC'\subset \cC$ is a morphism
of glutoses; the subglutos $(\cC',\cO')$ is a
{\bf full} subglutos of $(\cC,\cO)$ if $\cC'$
is a full subcategory of \cC.

\noindent {\bf Example:} Let \cE be a Grothendieck topos
(with respect to some universe \cU); let
$\cE^{-}$ be the full subcategory of \cE
consisting of all pointless objects, i.e.\ those
objects $X$ which have no global sections
$1 \arr X$. Then the pair
$(\cE^{-},\cE^{-})$ is a full subultraglutos
of the ultraglutos $(\cE,\cE)$. If one chooses
\cE properly, the glutos $(\cE^{-},\cE^{-})$
will have no terminal object (see example 0 of
section \ref{def}).

\noindent {\bf Counterexample:} If \cU is any universe containing some
infinite set, then the topos ${\bf Set}_{\cU_f}$
of finite sets is {\it not\/} a subglutos of the
topos ${\bf Set}_{\cU}$, because the corresponding
inclusion has no right adjoint.

\section{\cU-(ultra)glutoses}
\label{ug}
For any universe \cU there arise a
counterpart of Grothendieck toposes (=\cU-toposes by
terminology of \cite{SGA4}). Namely, call an ultraglutos $(\cC,\cO)$
an \cU-{\bf ultraglutos\/}, if \cC
is an \cU-category, any $\cO/X$ is an \cU-topos
and, besides, it satisfies the more strong axiom
{\bf(${\bf G4_{\cU}}$)}, obtained from the axiom {\bf(G4)} by replacing
`finite' by `\cU-small'.

To define an {\bf \cU-glutos} one needs as well to strengthen
the weak form of the axiom (G5), replacing in p.~c) of it
`finitely locally coequalizable'
by `locally coequalizable'.
(we will not distinguish in notations axiom~{\bf(G5)} for glutoses and
for \cU-glutoses).

By an \cU-{\bf category} is meant here a category with
\cU-small hom-sets which, besides, is naturally equivalent to some
category $\cC'$ with \mbox{Mor$\cC'\subset\cU$}
(i.e. this definition is stronger than that of \cite{SGA4}).

\begin{Rem}
One can show that p.b) of axiom
(G3) follows from other axioms for the case of
\cU-glutoses.
\end{Rem}

Examples (1)--(4) of section \ref{exa} above are examples of \cU-glutoses,
which are {\bf not} \cU-ultraglutoses;
another example is a functor category $\cC^{\cD}$,
where \cD is \cU-small and $(\cC,\cO)$ is
an \cU-glutos, if one defines the subcategory $\cO'$
in $\cC^{\cD}$ as follows:
\[\rho \colon F \arr F'\colon  \cD \arr\cC\]
belongs to $\cO'$ iff for any object $D$ of \cD
the arrow
\[\rho_D \colon FD \arr F'D \]
belongs to \cO. If \cC is an \cU-ultraglutos, then $\cC^{\cD}$ is
an \cU-ultraglutos as well.

Note that $\cO'$ is, generally speaking, bigger than $\cO^{\cD}$.




\begin{Rem}
\label{etstr}
\'Etale structures on an
\cU-topos \cE satisfying ``collection axiom'' in the sense of
\cite{JoMo},
are particular case of glutos structures on
toposes, as one can easily deduce from
Corollary 2.3 of \cite{JoMo}. Moreover, for any set
$Et$ of \'etale maps satisfying collection axioms
the pair $(\cE,Et)$ is a full subglutos
of the glutos $(\cE,\cE)$. As to interrelations
of glutos structures with \'etale structures, one can
see that any glutos structure \cO on an {\it arbitrary\/} category
\cC satisfies all of the conditions (A1)--(A8) of
\cite{JoMo} {\it excepting\/} conditions (A3) and (A6)
(one can easily find counterexamples in the glutos
{\bf Top}); and even in the case when \cC is an
\cU-topos
I have not discovered any special
relations (like ``descent'' and
``quotient'' axioms of \cite{JoMo}) between glutos structures \cO
on \cC and arbitrary epi's of \cC.
\end{Rem}


\section{Glueing in glutoses}
\label{glu}
In this section is studied what kind
of pullbacks and colimits exist in \cU-glutoses or\cU-ultraglutoses
(besides those whose existence is declared by axioms (G1), ${\rm(G4_{\cU})}$
and (G5)).

The general motto here is that in an \cU-glutos
pullback of two arrows exists if it exists locally
and that one can glue \cU-small families of objects
along open arrows. The rest of this section is devoted
to materialization of this motto into precise statements.

It turns out, in fact, that the corresponding
results are valid not only in \cU-glutoses (resp. in \cU-ultraglutoses),
but, more generally, in any clopos,
satisfying the corresponding version of axioms ${\rm(G4_{\cU})}$-(G5)
of section \ref{def} above.

For any set $I$ let $\Gamma I$ be the category defined
as follows. Its set of objects is just the set of all
non-empty words of length $\le2$ in the free monoid
$W(I)$ of $I$-words (which is supposed to be chosen in such
a way that $I$ is a subset of $W(I)$ and the canonical
map $I \arr W(I)$ coincides with the inclusion of subsets). The only non-identity arrows of $\Gamma I$ are
arrows
\[ i \stackrel{d_0^{ij}}{\longleftarrow} ij
\stackrel{d_1^{ij}}{\arr} j\ \ (i,j\in I); \]
note that $d_0^{ii}$ is to be {\it different} from
$d_1^{ii}$.

Given a diagram $U\colon \Gamma I \arr \cC$
we will write $U_i$, resp. $U_{ij}$, resp. $d_\varepsilon^{ij}$
instead of $U(i)$, resp. $U(ij)$, resp. $U(d_\varepsilon^{ij})$
omitting sometimes superscripts in the latter case;
if $U$ has a colimit we will write \mbox{$U.$} for a
colimit object and $\{ \mbox{$u_{i}.$}\colon U_i \arr
\mbox{$U.$} \}_{i\in I}$ for a colimit cone.

Call a diagram $U\colon \Gamma I \arr \cC$
{\bf glueing data} or a {\bf gluon} if the following ``cocicle conditions'' are satisfied:

\noindent {\bf(GD1)} For any $i,j\in I$ the pair $(d_0^{ij},d_1^{ij})$ is a
mono source in \cC;

\noindent {\bf(GD2)} For any $i,j\in I$ there exists
an arrow $\tau_{ij}\colon U_{ij} \arr U_{ji} $
of \cC such that the equalities
\[ d_0^{ji}\tau_{ij} = d_1^{ij} \ \ {\rm and} \ \
d_1^{ji}\tau_{ij} = d_0^{ij};\]
are valid (it then follows from (GD1) that
$\tau_{ij} \tau_{ji} = {\rm Id}$);

\noindent {\bf(GD3)} For any $i\in I$ there exists
an arrow $s_i \colon U_i \arr U_{ii}$
of \cC which is right inverse to both
$d_0^{ii}$ and $d_1^{ii}$;

\noindent {\bf(GD4)} For any word $ijk$ of length 3
in $W(I)$ there exists an object $U_{ijk}$ of \cC
and the arrows
$p_0\colon U_{ijk} \arr U_{ij}$,
$p\colon U_{ijk} \arr U_{ik}$ and
$p_1\colon U_{ijk} \arr U_{jk}$
such that all three squares of the
diagram
\begin{equation}
\label{pic1}
\expandafter\ifx\csname xymatrix\endcsname\relax%
\input{pic1}
\else
\xymatrix{
&U_{ijk}\ar[dl]_{p_0}\ar[d]^p\ar[dr]^{p_1}\\
U_{ij}\ar[d]_(.45){d_0}\ar[dr]_(.3){d_1}&
U_{ik}\ar[dl]^(.3){d_0}|\hole\ar[dr]_(.3){d_1}|\hole&
U_{jk}\ar[dl]^(.3){d_0}\ar[d]^(.45){d_1}\\
U_i&U_j&U_k
}
\fi
\end{equation}
are pullbacks (we will write further $p_0^{ijk}$, etc. instead of $p_0$ in case of necessity).

It follows from the latter condition the existence of
isomorphisms $\theta_{ijk}\colon U_{ijk} \arr
U_{jki}$ which agree with projections $p_0$, $p_1$
and $p$ ``twisted'' by iso's $\tau_{ij}$ and satisfy
the ``cocicle conditions'' arising both in algebraic and
differential geometry in processes of glueing of
schemes, resp. manifolds along open subschemes, resp.
open submanifolds.

One can see, on the other hand, that if the index set $I$
consists of just one element the definition above
reproduces the definition of an equivalence relation.

More generally, for any glueing data $U\colon \Gamma I \arr \cC$
and any $i\in I$ the pair \darr{d^{ii}_0}{d^{ii}_1}{U_{ii}}{U_i}
is, evidently, an equivalence relation.

Any family $\{u_i\colon U_i \arr X\}_{i\in I}$ of pullbackable
arrows defines, canonically, some glueing functor if one sets
$U_{ij}\approx U_i \prod_X U_j$, whereas for $d^{ij}_\varepsilon$
one chooses the corresponding pullback projections.

Call a diagram in a clopos (resp. in a glutos) {\bf clopen}
(resp. {\bf open}) if any arrow
of this diagram is clopen (resp. open).
Call a clopen gluon $U\colon \Gamma I \arr \cC$ in a clopos $(\cC,\cO)$
{\bf(finitely) locally \cO-coequalizable} if ($I$ is finite and)
for any $i\in I$ the equivalence relation
\darr{d^{ii}_0}{d^{ii}_1}{U_{ii}}{U_i}
is (finitely) locally \cO-coequalizable.

One sees immediately that
for any clopen gluon in a clopos, morphisms $\tau_{ij}$,
$p_0$, $p_1$ and $p$ of (GD2), (GD4) are clopen. As to
(the only by (GD1)) arrows $s_i$ of (GD3), they also
are clopen if the clopos satisfies conditions ${\rm(G4_{\cU})}$-(G5)
as one can see from the following

\begin{Prop}
\label{glufun}
Suppose that a clopos $(\cC,\cO)$ satisfies conditions ${\rm(G4_{\cU})}$-{\rm(G5)}
in the definition of \cU-ultraglutoses {\rm(resp.} of \cU-glutoses\/{\rm).}
Then:

\noindent {\bf(${\bf G4_{\cU}+5}$)} Any \cU-small clopen gluon
{\rm(resp.} any \cU-small clopen locally \cO-coequalizable gluon\/{\rm)}%
\footnote{Note that any clopen gluon every arrow of which is mono,
is locally \cO-coequalizable automatically, because corresponding
equivalence relations are trivial in this case:
both $d_0$ and $d_1$ are iso's.}
$U\colon \Gamma I \arr \cC$
has a universal colimit $U.$
which, besides, is effective in the sense
that for any $i,j\in I$ the isomorphism
\[U_{ij} \approx U_i \prod_{U.} U_j \]
holds. The colimit cone
$\{U_i \arr U.\}_{i\in I}$
consists of clopen arrows.
\end{Prop}

\noindent{\bf Indications to the proof.} Consider the diagram
\begin{equation}
\label{eqreldia}
d_0,d_1\colon \coprod_{i,j\in I}U_{ij}
{}^{\arr}_{\arr} \coprod_{i\in I}U_i,
\end{equation}
where, say, $d_0$ is defined as $(\iota_id^{ij}_0)_{i,j\in I}$ with
$\iota_i\colon U_i \arr \coprod_{i\in I}U_i$ being the
canonical coproduct injection arrows\footnote{We use parentheses
instead of braces in order to distinguish between
{\it families} of arrows and a single (co)limit arrow
determined by the corresponding family.}.
One is to prove that the
diagram above is an open equivalence relation
(resp. is an open \cO-coequalizable equivalence relation);
then it will follow trivially that the coequalizer
$$q\colon \coprod_{i\in I}U_i \arr \mbox{$U.$}$$
of this diagram (existing by ${\bf(G5_{\cU})}$)
reproduces the colimit cone of the original
gluon $U$ if one sets $\mbox{$u_i.$}=q\iota_i$.

Note first that the families $\{\tau_{ij}\}_{i,j\in I}$ and
$\{s_i\}_{i\in I}$ of (GD2) and (GD3) permit one to build
in a natural way the arrows $\tau\colon
\coprod U_{ij}\arr
\coprod U_{ij}$
and
$s\colon \coprod U_i \arr
\coprod U_{ij}$; the verification
of the fact that these arrows satisfy (GD2), resp.(GD3),
is straightforward.

Similarly, one can construct three arrows $p_0,\ p_1$ and $p$ from $\coprod U_{ijk}$ to
$\coprod U_{ij}$; e.g., the arrow $p_0\colon
\coprod U_{ijk}\arr \coprod U_{ij}$ is defined to be the colimit arrow
\[(U_{ijk}\stackrel{p_0^{ijk}}{\arr}U_{ij}
\stackrel{\iota_{ij}}{\arr}\coprod U_{ij})_{i,j,k\in I}\ \ ,\]
where $\iota_{ij}$ are canonical coproduct injection arrows.

In proving  (GD4) for the diagram~\ref{eqreldia} above the following useful lemma can be used, which
states that a square is a pullback if it is a pullback
locally:

\begin{Lem}
\label{pullcriter}
Let $f\colon X \arr Z$ and $g\colon Y \arr Z$ be arrows of a category \cC;
let $\{x_i\colon X_i \arr X\}_{i\in I}$,
$\{y_j\colon Y_j \arr Y\}_{j\in J}$ and
$\{z_k\colon Z_k \arr Z\}_{k\in K}$ are
universal effective epi families; let, further, for
any $i\in I$, $j\in J$ and $k\in K$ a diagram
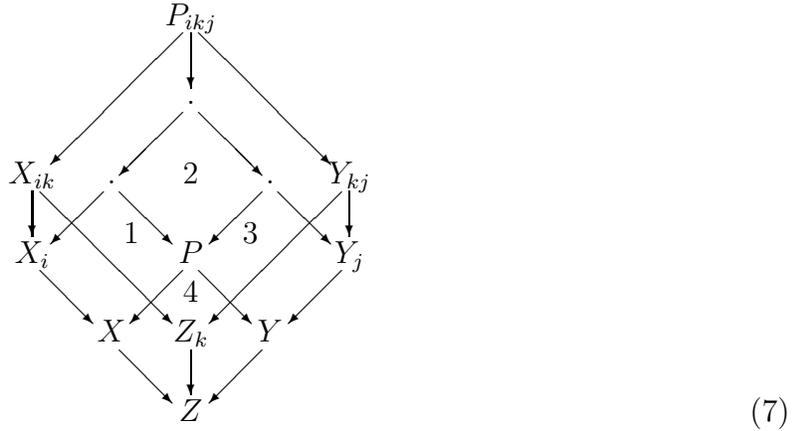
\def\ar#1 #2,#3 #4,#5;{\put(#1,#2){\vector(#3,#4){#5}}}
\def\bx#1 #2,#3;{\put(#1,#2){\makebox[0pt]{$#3$}}}
\def\ln#1 #2,#3 #4,#5;{\put(#1,#2){\line(#3,#4){#5}}}
\def\mar#1 #2,#3 #4,#5;#6 #7,#8;{\multiput(#1,#2)%
(#6,#7){#8}{\vector(#3,#4){#5}}}
\def\mln#1 #2,#3 #4,#5;#6 #7,#8;{\multiput(#1,#2)%
(#6,#7){#8}{\line(#3,#4){#5}}}
\begin{equation}
\label{pic2}
\begin{picture}(120,160)
  \bx60 0,Z;
  \bx30 30,X;\bx60 30,Z_k;\bx90 30,Y;
  \bx0 60,X_i;\bx60 60,P;\bx120 60,Y_j;
  \bx0 90,X_{ik};\bx30 90,.;\bx90 90,.;\bx120 90,Y_{kj};
  \bx60 120,.;
  \bx60 150,P_{ikj};
  \mar0 87,0 -1,17;60 -60,2;
  \ar60 147,0 -1,21;\ar120 87,0 -1,17;
  \mar3 57,1 -1,20;30 30,2;\ar63 117,1 -1,24;
  \mar33 27,1 -1,20;30 30,3;
  \ar57 117,-1 -1,24;\mar87 87,-1 -1,20;30 -30,2;
  \mar27 87,-1 -1,20;30 -30,3;
  \mar3 87,1 -1,50;60 60,2;
  \mar57 147,-1 -1,50;60 -60,2;
  \bx37.5 67.5,1;\bx60 90,2;\bx82.5 67.5,3;\bx60 45,4;
\end{picture}
\end{equation}
is given such that four side squares of it as well as
the ``floor'' squares 1,~2 and~3 are pullbacks. Then
the square 4 is a pullback iff for any $i\in I$, $j\in J$ and $k\in K$ the ``ceiling'' square is a pullback.
\end{Lem}

Applying this lemma to the case $X=Y=\coprod U_{ij}$,
$Z=\coprod U_i$ and $P=\coprod U_{ijk}$ with the
corresponding universal effective epi families being
$\{\iota_{ij}\colon U_{ij}\arr \coprod U_{ij}\}_{i,j\in I}$, etc., one obtains after simple diagram chase just squares
of the diagram~(\ref{pic1}) as  ``ceilings'' of the diagram~(\ref{pic2})
above, which proves (GD4) for the diagram~(\ref{eqreldia}).

At last, to prove (GD1) for the diagram~(\ref{eqreldia})
consider a pair of arrows
$f,g\colon X{}_{\arr}^{\arr} \coprod U_{ij}$ such that both
$d_0f=d_0g$ and $d_1f=d_1g$. Pulling a covering
$\{\iota_i\colon U_i \arr \coprod U_i\}_{i\in I}$
along $d_0f=d_0g$, resp. along $d_1f=d_1g$, one
obtains two universal effective epi families
$\{v_i\colon V_i \arr X\}_{i\in I}$ and
$\{v'_j\colon V_j \arr X\}_{j\in I}$ such that $f$ agree
with $g$ on elements of the ``intersection'' universal effective epi family
$\{V_i\coprod_X V'_j \arr X\}_{i,j\in I}$, which implies that
$f=g$.\
\eop

The following proposition describes sufficient conditions
of existence of pullbacks in cloposes satisfying ${\rm(G4_{\cU})}$-(G5).

\begin{Prop}
\label{pull}
Let a clopos $(\cC,\cO)$ satisfies conditions
${\rm(G4_{\cU})}$ as well as the strong axiom {\rm(G5)} or
the combination of the weak axiom~{\rm(G5)} with the axiom~{\rm(G5P)}.
Let $f\colon X \arr Z$ and
$g\colon Y \arr Z$ be arrows of \cC such that
for some \cU-small epi families of clopen arrows $\{X_i \arr X\}_{i\in I}$, $\{Y_j \arr Y\}_{j\in J}$ and
$\{Z_k \arr Z\}_{k\in K}$ there exists, for any
$i\in I$, $j\in J$ and $k\in K$, the pullback
$X_{ik}\prod_{Z_k} Y_{kj}$, where, by definition,
$X_{ik}=X_i\prod_Z Z_k$ and
$Y_{kj}=Z_k\prod_Z Y_j$. Then there exists the pullback
of $f$ and $g$. If, besides, $(\cC,\cO)$ is an
\cU-glutos, then the \cU-smallness condition
for families above can be omitted.
\end{Prop}

The archetype of the proof of Prop. \ref{pull} is contained,
for example, in the proof of existence of pullbacks of
Grothendieck schemes (see, e.g., \cite{Har}).

Proposition \ref{glufun} permits one to equip, canonically, any \cU-glutos $(\cC,\cO)$ with a structure of a presite, but, before going into details, one needs to give
some necessary definitions and to state some elementary properties of presites.

\section{Presites}
\label{ps}
Define first, for a presite
$(\cC,\tau)$ the set $\cO_{\tau}$ of arrows of \cC
as consisting of just those arrows
which belong to some covering of $\tau$. The set $\cO_{\tau}$
satisfies condition (G1) above (so that
its elements will be referred to as clopen or open) and
one can define {\bf morphisms} between presites
$(\cC,\tau)$ and $(\cC',\tau')$ as just those
functors between underlying categories which respect coverings and satisfy
condition (MG1) above (with \cO replaced by $\cO_{\tau}$,
idem for $\cO'$). We will call
such functors {\bf continuous} (this definition is
stronger than the corresponding definition in \cite{SGA4}
making emphasis on topologies and sites).

If $F\colon\cC \arr\cC'$ is a functor
and $\tau'$ is a pretopology on $\cC'$, then
a pretopology $\tau$ on \cC will be called
{\bf induced} by $\tau'$ along $F$ iff for any sink $S$
in \cC the condition $FS\in \tau'$ is equivalent to $S\in \tau$;
if such $\tau$ exists it is the biggest
pretopology on \cC making $F$ continuous.

\begin{Prop}
\label{ind}
For any presite $(\cC,\tau)$ and any object
$X$ of \cC there exists the pretopology on $\cO_{\tau}/X$ induced by
$\tau$ along the ``source'' functor
$d_{0}\colon\cO_{\tau}/X \arr\cC$.
\end{Prop}

The category $\cO_{\tau}/X$ will be considered, canonocally,
to be equipped with this presite structure; then:

\begin{Prop}
\label{indfun}
For any arrow $f\colon X \arr Y$ of \cC the induced
functor $f^{*}\colon\cO_{\tau}/Y \arr\cO_{\tau}/X$
is continuous.
\end{Prop}

Now, a presite $(\cC,\tau)$ will be called
an \cU-{\bf presite} if \cC is an
\cU-category and, besides, the following
condition is satisfied (existence of local sets of
topological generators):

\noindent {\bf(${\bf LTG_{\cU}}$)} For any object $X$ of \cC
there exists an \cU-small subset $G_{X}$ of
objects of \cC such that for any clopen arrow
$u\colon U \arr X$ of \cC there exists
a covering \mbox{$\{ u_i\colon U_i \arr U\}\iI$}
such that any $U_i$ belongs to $G_X$. This condition is just equivalent
to saying that any $\cO_{\tau}/X$, considered
as a site, is an \cU-site in the sense of \cite{SGA4}.

\begin{Rem}
\label{pt4}
It is convenient to include in the definition of a pretopology the following
condition (completeness property):

\noindent {\bf(PT4)} If $\langle X,S\rangle$ is a sink of \cC such that
$S\subset \cO_{\tau}$ (such sinks will be called {\bf (cl)open})
and there exists a refinement of $\langle X,S\rangle$ which is a
covering of $X$ then $\langle X,S\rangle$ itself is a covering. Here
a sink $\langle X,S')$ is said to be a {\bf  refinement} of
$\langle X,S\rangle$
if for any $s'\in S'$ there exists $s\in S$ such that $s'$
factors through $s$.

Intersection of pretopologies satisfying (PT4) satisfies
(PT4) itself; if $\tau$ satisfies (PT4) then a pretopology
induced by $\tau$ along any functor satisfies (PT4) as well;
besides, the completion of a pretopology
$\tau$ satisfying ordinary
conditions (PT1)--(PT3) of \cite{SGA4} to that satisfying
(PT4), does not change neither  the set
$\cO_{\tau}$, neither the associated Grothendieck
topology%
\footnote{This completion is, in fact, the biggest pretopology among
those having both the same set of open arrows and the same associated
\grot topology as $\tau$ has.},%
nor the universal glutos of Theorem~\ref{th} below.
That is why from now on `pretopology' will
mean `pretopology satisfying (PT4)' with similar
change in the meaning of `presite'.
\end{Rem}

\begin{Rem}
If one looks at the definition of elementary glutos, a natural
question can arise: what will happen if one ``iterates''
the theory of glutoses replacing, roughly, in axioms
(G1)--(G5) ``topos'' by ``glutos''? The answer is that
the theory of elementary glutoses is stable by this
iteration, i.e., no new ``weaker'' theory will arise.

In more details, defining, in an evident way, morphisms
of cloposes as well as clopos structure induced along
a functor, one obtains that for any object $X$ of any
clopos $(\cC,\cO)$ the category $\cO/X$ has
a clopos structure $\cO_X$ induced along the functor
$d_0$: arrows of $\cO_X$ are all commutative triangles
(i.e., arrows of $\cO/X$) all three arrows of which
belong to \cO. Counterparts of Props. \ref{ind} and~\ref{indfun}
are valid for cloposes as well as the following result:

\begin{Prop}
\label{iter}
For any clopen arrow $U\stackrel{u}{\arr}X$
in a clopos $(\cC,\cO)$ the functor
\[
d_0/(U\stackrel{u}{\arr}X)\colon
\cO_X/(U\stackrel{u}{\arr}X)\arr
\cO/U
\]
is a natural equivalence.
\end{Prop}

Now, if one {\it removes\/} the axiom (G2) and one replaces
in axiom (G3) ``topos'' by ``glutos'', resp. ``inverse image
of geometric morphism'' by ``morphism of glutoses'', interpreting,
simultaneously, $\cO/X$, etc. not simply as categories
but as cloposes via induced structure, then one arrives
to an elementary theory which turns out to be not weaker, but
equivalent
to the theory of elementary glutoses. This just follows
from Prop. \ref{iter}.
\end{Rem}

\section{\cU-glutoses as \cU-presites}
\label{ugps}
Returning again to \cU-glutoses, one has:

\begin{Prop}
\label{proper}
Let $(\cC,\cO)$ be an \cU-glutos. Then:

\noindent {\bf(${\bf G6_{\cU}}$)} All epi sinks in \cC with elements
in \cO are universal effective and, hence, form some pretopology
on \cC {\rm (denoted further $\tau_{\cO}$).
This pretopology is {\bf subcanonical\/}
(i.e. the associated topology is subcanonical);}

\noindent {\bf(${\bf G7_{\cU}}$)} The presite $(\cC,
\tau_{ \cO})$ is an \cU-presite;

\noindent {\bf(${\bf G8_{\cU}}$)} Any sink $\langle X,S\rangle$ with
$S \subset \cO$ factors as a covering of
$\tau_{\cO}$ followed by an open mono.

\noindent {\bf(${\bf G9_{\cU}}$)} {\rm(local character of open arrows)}
If for $u\colon U \arr X$ there exists
a covering $\{ u_i\colon U_i \arr U \}\iI$
of $\tau_{\cO}$
such that for any $i\in I$ the arrow $uu_i$ is open,
then $u$ itself is open.

\noindent {\bf(G10)} For any object $X$ of \cC
the pretopology on $\cC/X$ induced by $\tau_{\cO}$
is the canonical pretopology of the topos $\cC/X$
{\rm (i.e. coverings of it are all epi sinks);}
moreover, the source functor
$d_0\colon  \cC/X \arr\cC$
respects both coequalizers of equivalence relations
and \cU-small coproducts.
\end{Prop}

\begin{Rem}
For an elementary glutos
$(\cC,\cO)$ let $\tau_{\cO}$ be the set
of all open sinks having some finite open epi refinement.
Then one has: {\bf(G6)} $\tau_{\cO}$ is a subcanonical
pretopology on \cC; the counterparts of ${\rm(G8_{\cU})}$
and ${\rm(G9_{\cU})}$ are valid as well if one replaces
in ${\rm(G8_{\cU})}$ `Any sink' by `Any finite sink'.
\end{Rem}

The following proposition is a counterpart of Giraud
theorem :

\begin{Prop}
\label{gir}
A pair $(\cC,\cO)$ is an \cU-glutos {\rm(resp.} an \cU-ultraglutos\/{\rm)}
iff it satisfies conditions
{\rm (G1)--(G2), ${\rm(G4_{\cU})}$, {\it weak axiom\/} (G5) (resp.,
{\it strong axiom\/} (G5) {\it and axiom\/} (G5P)),
${\rm(G6_{\cU})}$-${\rm(G7_{\cU})}$} above
{\rm(conditions ${\rm(G4_{\cU})}$ (G5) and (G5P) can be
replaced by weak (resp. strong) form of condition ${\rm(G4_{\cU}+5)}$ and
condition ${\rm(G9_{\cU})}$).}
\end{Prop}

Now, a map $(\cC,\cO) \mapsto (\cC,\tau_{\cO})$
continues to the 2-functor imbedding fully \cU-glutoses
into \cU-presites, as shows the following

\begin{Prop}
\label{cont}
Let $(\cC,\cO)$
and $(\cC',\cO')$ be \cU-glutoses and
$F\colon\cC \arr\cC'$ be a functor. Then
$F$ is morphism of glutoses iff it is continuous
w.r.t. pretopologies $\tau_{\cO}$ and
$\tau_{\cO'}$.
\end{Prop}

In other words, the
2-category \mbox{{\bf Glut}${}_{\cU}$} of \cU-glutoses
may be considered as a full
2-subcategory of the 2-category
\mbox{{\bf Psite}${}_{\cU}$} of \cU-presites
\begin{equation}
\label{glusubps}
\mbox{{\bf Glut}${}_{\cU}$} \harr
\mbox{{\bf Psite}${}_{\cU}$}.
\end{equation}

\begin{Rem}
Let \mbox{{\bf Glut.}},
be the full 2-subcategory
of {\bf Glut},
containing any glutos which is \cU-glutos for some
universe \cU. Let the 2-category \mbox{{\bf Psite.}}
be defined similarly. The above inclusion functor continues to
the inclusion functor
\[ \mbox{\bf Glut.} \harr\mbox{\bf Psite.}\ , \]
but the counterexample
\[{\bf Set}_{\cU_f} \harr {\bf Set}_{\cU} \]
of continuous functor which is not a morphism of glutoses
(see the end of section \ref{2cat})
shows that this inclusion is {\it not\/} full.
\end{Rem}

The main author's result states that {\it the 2-category
\mbox{{\bf Glut}${}_{\cU}$} of \cU-glutoses is reflective in the
2-category \mbox{{\bf Psite}${}_{\cU}$} of \cU-presites, whereas
the 2-category \mbox{{\bf UGlut}${}_{\cU}$} of \cU-ultraglutoses
is reflective in the 2-category \mbox{{\bf Glut}${}_{\cU}$}.}
In more details:
\begin{Th}
\label{th}
{\rm(a)} For any \cU-presite
\cC there exists an \cU-(ultra)glutos $\wt{\cC}$
to\-ge\-ther with an ar\-r\-ow
\mbox{$Y_{\cC}\colon\cC\arr\wt{\cC}$}
which is universal in the sense that for any \cU-glutos \cD the arrow
\begin{equation}
\label{natequiv}
[Y_{\cC},\cD]\colon [\wt{\cC},\cD]
\arr [\cC,\cD] \ \ (\sigma \mapsto
\sigma Y_{\cC})
\end{equation}
is a natural equivalence having right inverse;

\noindent {\rm(b)} the arrow $Y_{\cC}$ {\rm (or, rather,
the underlying functor)} can always be chosen to be injective on objects
of \cC;
if the pretopology $\tau$ of \cC is subcanonical,
then the arrow $Y_{\cC}$ is fully faithful;

\noindent {\rm (c)} The universal arrow $Y_{\cC}$
reflects open coverings {\rm (see sect. \ref{charts}
below for the definition)}; for every object
$X$ of $\wt{\cC}$ the set $G_X$ of topological
generators of $X$ {\rm (see ${\rm(LTG_{\cU})}$ in sect. \ref{ps}
above)} can be chosen to belong to the set $Y_{\cC}(\cC)$;

\noindent {\rm (d)} Besides, if the underlying category
of \cC {\rm (denoted further \cC, by abuse of notation)}
is \cU-cocomplete and the pretopology of \cC
is subcanonical, then the functor $Y_{\cC}$
has left adjoint
\mbox{$\Gamma \colon  \wt{\cC} \arr\cC$}
{\rm (the {\bf global sections functor}). Note that
$\Gamma$ need not be
continuous.}
%
\end{Th}

The proof of Th.\ref{th} is sketched in
Appendix~\ref{Proof}.

\begin{Rem}
It follows from
Appendix~\ref{Proof}
that the 2-category
of subcanonical \cU-presites is as well reflective in
\mbox{{\bf Psite}${}_{\cU}$}, i.e. the universal arrow
$Y_{\cC}$ decomposes as
\[\cC\stackrel{Y'_{\cC}}{\arr}\cC_{sub}
\stackrel{Y_{\cC_{sub}}}{\arr}\wt{\cC_{sub}},\]
where $\cC_{sub}$ is a universal subcanonical \cU-presite for \cC.
\end{Rem}


Now choosing for every \cU-presite \cC some universal
arrow $Y_{\cC}$ and choosing for every pair
\cC,~\cD as in p.(a) of Theorem \ref{th} some arrow
\[ I_{\cC\cD}\colon [\cC,\cD] \arr
[\wt{\cC},\cD] \]
right inverse to the arrow (\ref{natequiv}),
we will obtain for every pair \cC,~$\cC'$ of
\cU-presites some functor
\[ [\cC,\cC']\stackrel{\sim}{\arr}
[\wt{\cC},\wt{\cC'}]\ \
(\sigma\mapsto\wt\sigma)\ ,\]
defined by $\wt\sigma:=I_{\cC\wt{\cC'}}
(Y_{\cC'}\sigma)$ on 2-arrows from $[\cC,\cC']$.

The correspondences $\sigma\mapsto\wt\sigma$ just
defined are incorporating together to give some pseudofunctor
(see \cite{SGA4},\cite{Gra})
\[\mbox{{\bf Psite}${}_{\cU}$}
\stackrel{\sim}{\arr}
\mbox{{\bf Glut}${}_{\cU}$}\ ,\]
left quasiadjoint to the inclusion 2-functor~(\ref{glusubps});
it differs from a 2-functor by some ``twisting by a cocicle''
$\sigma(F,F')\colon\wt{F'}\wt{F}
\arr\wt{(F'F)}$. The following theorem
shows that this cocicle can be killed.

\begin{Th}
\label{th2}
The correspondences $\cC\mapsto\wt{\cC}$ and
$F\mapsto\wt{F}$ can be chosen in such a way that
$\wt{F'}\wt{F}=\wt{(F'F)}$ for every
composable pair $F$~and~$F'$ of morphisms of
\cU-presites.
\end{Th}

\begin{Cor}
If $F\colon\cC\arr\cC'$ and
$G\colon\cC'\arr\cC$ are continuous
functors between \cU-presites such that $F$ is left
adjoint to $G$ then $\wt{F}$ is left adjoint to
$\wt{G}$.
\end{Cor}

Exactness properties of universal arrows are described by
the following

\begin{Prop}
\label{exact}
{\rm (a)} For any \cU-presite \cC
the universal arrow $Y_{\cC}\colon\cC\arr
\wt{\cC}$ respects all \cU-li\-m\-its e\-xi\-s\-t\-ing
in \cC;

\noindent {\rm (b)} If \cC has pullbacks, resp.
products, resp. finite limits, then so does
$\wt{\cC}$;

\noindent {\rm (c)} Let $F\colon\cC \arr\cC'$
be morphism of presites and \cC has products, resp.
pullbacks, resp. finite limits, which are, besides, respected by the functor
$F$. Then the functor
$\wt{F}\colon  \wt{\cC} \arr\wt{\cC}'$
respects products, resp. pullbacks, resp. finite limits.
\end{Prop}

Theorem \ref{th} and Prop.\ref{exact} show that glutoses are
``invariants'' of presites in just the same way as toposes
are ``invariants'' of sites.
The universal arrow for an \cU-presite is, clearly,
a counterpart of topos-theoretic ``sheafified Yoneda functor''
$Y\colon\cS \arr {\rm Sh}\cS$ associating
to any \cU-site \cS the topos of ${\bf Set}_{\cU}$-valued sheaves
on it.
In many familiar cases of \cU-presites \cC (see examples
of universal arrows below), the corresponding {\em site}
is not an \cU-site, which means that the topos of sheaves
on this site exists in some {\em higher} universe only. At the same time,
the glutos $\wt{\cC}$ associated with the presite \cC
exists in the {\em same} universe \cU, where \cC is contained.
Nevertheless, when both constructions exist, they
sometimes coincide as shows the following

\begin{Prop}
\label{sh}
Let \cC be
\cU-small and finitely complete. Let a pretopology
$\tau$ on \cC be given, such that any arrow of \cC
is clopen. Then the universal glutos $\wt{\cC}$
coincides with the topos of sheaves {\rm Sh}\cC
up to natural equivalence of categories. The similar is true for
universal arrows.
\end{Prop}

For example, the glutos constructed from a topological
space is the same thing as the topos of sheaves on it;
the same is true for a complete Heyting algebra
(equipped with the canonical (pre-)topology).

Many familiar examples of constructing categories out of
``simpler ones'' by means of ``charts and atlases'' routine
are just concrete realizations of universal arrows of Theorem \ref{th}:
imbedding of smooth euclidean regions into
the category of smooth manifolds, imbedding of trivial
vector bundles into the category of locally trivial ones,
as well as the functor
\[ {\rm Spec\colon }{\bf Ring}^{op} \arr {\bf Schem}.\]

Note that in this example the global sections functor of
p.(d) of Theorem \ref{th} exists and is the same thing as
the ordinary global sections functor on {\bf Schem}, which
justifies the use of the name ``global sections functor''
in the general case.

The latter example opens up a new approach to ``universal
algebraic geometry'', alternative to that of M.Coste
\cite{Cos} (based on Hakim's theorem): given some
locally finitely presentable category (see \cite{GaUl}) \cC
together with some pretopology $\tau$ on its
dual category, turning $\cC^{op}$ into an
\cU-presite, the category of {\bf schemes over} \cC
and the corresponding functor Spec can be {\em defined}
to be, respectively, the glutos associated with the
presite $(\cC^{op},\tau)$ and the universal
arrow for it.

For example, if one chooses the \'{e}tale pretopology
on the category dual to that of commutative rings
instead of Zariski pretopology, one obtains the category
{\bf ESchem}, which may be called the category of
{\bf \'{e}tale schemes}; the functor
\begin{equation}
\label{scharresch}
{\bf Schem} \arr {\bf ESchem},
\end{equation}
provided by Theorem \ref{th}, fully imbeds the category
of schemes into that of \'{e}tale schemes in such a way
that for any scheme $X$ the topos of sheaves over $X$
with respect to \'etale pretopology on $Et/X$ imbeds
into {\bf ESchem} via
\[
{\rm Sh}(Et/X) \harr {\bf ESchem}/X
\stackrel{d_0}{\arr} {\bf ESchem}.
\]

Another application is concerned with ``non-commutative
algebraic geometry'': Theorem \ref{th} gives general
non-commutative schemes glued out of non-commutative
affine schemes of P.M.Cohn \cite{Coh}.

A class of pretopologies on duals to locally finitely
presentable categories especially suitable for ``universal
algebraic geometry'' will be considered elsewhere.

\section{\cM-presites and SG-glutoses.}
\label{oper}
In this section
some natural endo-2-functors are constructed on the
2-category {\bf Psite} of all presites; recall,that
\cU-glutoses are considered as presites via Props.
\ref{proper}~and~\ref{cont}.
Besides, a class of \cU-glutoses  which
occur particularly often in practice is studied in more details.
%

Let \cP be some property of an arrow of a presite \cC.
We say that an arrow $f\colon X \arr Y$ of \cC
{\bf locally satisfies} \cP or {\bf is locally} \cP,
if there exists a covering $\{u_i\colon U_i \arr
X\}_{i\in I}$ such that for every $i\in I$ the arrow $fu_i$
satisfies \cP. We will use further this metadefinition
for the case when the property \cP is either ``$f$ is
(cl)open'' or ``$f$ is an (cl)open mono'' getting the properties
``$f$ is locally (cl)open'' or ``$f$ is locally an (cl)open mono''
(note yet that an $f$ which locally is a clopen mono
need not to be neither clopen nor mono). One can easily
verify that the set of all {\it pullbackable\/} locally open arrows
of any presite is closed both with respect to compositions and
arbitrary pullbacks.

For example, the property ${\rm(G9_{\cU})}$ of glutoses
(see section \ref{ugps}) can be reformulated in this terms as
follows: every locally open arrow in a glutos is open.

Let $\tau$ be a pretopology on a category \cC.
Define the pretopology $\cM\tau$, resp.
$\cL\tau$, resp. ${\rm SG}(\tau)$ on the category \cC
as follows: coverings of $\cM\tau$ are all coverings
of $\tau$ consisting of mono's; coverings of $\cL\tau$
are all sinks of $\tau$ consisting of
pullbackable locally open arrows of \cC and having
a refinement belonging to a pretopology $\tau$ (the
latter definition is correct because pullbackable
locally open arrows form a clopos structure as stated above);
at last let ${\rm SG}(\tau)=(\cL\cM\tau)\cap\tau$.
One has, evidently, the following inclusions:
\[\cM\tau\subset{\rm SG}(\tau)\subset\tau\subset
\cL\tau\ .\]

The operations $\cM$,~$\cL$ and~${\rm SG}$
can be continued to the endo-2-functors (denoted by
the same symbols) on the 2-category {\bf Psite},
whereas the chain of inclusions above produce the chain
of 2-functor morphisms
\begin{equation}
\label{fmm}
\cM\harr{\rm SG}\harr{\rm Id}_{\bf Psite}\harr\cL\ ,
\end{equation}
which go to identity 2-functor morphisms when being composed with
the neglecting 2-functor from {\bf Psite} to the 2-category
{\bf Cat} of all categories.

It is evident that for any \cU-presite \cC the presite
$\cL\cC$ is an \cU-presite (but the author do not know
at present whether or not the same is true for $\cM\cC$
and~${\rm SG}\cC$). An evident fact that the Grothendieck
{\it topologies\/} generated by pretopologies of \cC and of
$\cL\cC$ coincide, imply, together with the construction
of universal glutoses out of the corresponding category
of ``big'' sheaves
(see
Appendix~\ref{Proof}),
the following

\begin{Prop}
For any \cU-presite \cC there is a canonical natural equivalence
$\widetilde{\cC}\approx\widetilde{\cL\cC}$;
in more details, the composition arrow
\[\cC\harr\cL\cC
\stackrel{Y_{\cL\cC}}{\arr}\widetilde{\cL\cC}\]
is a universal arrow for \cC.
\end{Prop}

The following proposition describes the monoid of endo-2-functors
generated by $\cM$,~SG and~$\cL$.

\begin{Prop}
\label{rels}
The {\rm2}-functors $\cM$, {\rm SG} and $\cL$ satisfy the
following algebraic relations:
\begin{eqnarray*}
{\cM^2=\cM,}&{\rm SG}^2={\rm SG},&{\cL^2=\cL,}\\
{(\cM\cL)^2=\cM\cL,}&{(\cL\cM)^2=\cL\cM,}&{{\rm SG}\cL=\cL\cM\cL,}\\
{\cL{\rm SG}=\cL\cM,}&{{\rm SG}\cM=\cM,}&
{\cM{\rm SG}=\cM.}
\end{eqnarray*}
\end{Prop}

The only relation amongst those above, whose verification
uses drawing of some
diagrams is that stating the idempotence of the functor \cL.

The first three relations of proposition \ref{rels} together
with the universality properties of functor morphisms (\ref{fmm})
imply that both the full 2-subcategory of presites stable by \cM and
of presites stable by SG are coreflective in {\bf Psite},
whereas the full 2-subcategory of presites stable by \cL
is reflective in {\bf Psite}.

A presite stable by $\cM$, resp, by {\rm SG}, resp. by $\cL$ will be
called an \cM{\bf-presite}, resp. an {\rm SG}{\bf-presite},
resp. an \cL{\bf-presite}.
In other words, a presite \cC is an \cM-presite
iff any covering of it consists of mono's; it is an
{\rm SG}-presite iff any clopen arrow of it is locally
a clopen mono;
it is an \cL-presite iff any arrow
of it which is locally clopen is clopen.

In particular, any \cU-glutos is an \cL-presite;
an \cU-glutos $(\cC,\cO)$ is an SG-presite iff for any object
$X$ of \cC the topos $\cO/X$ is an SG-topos
as defined in \cite{Joh}, which justifies the name ``SG-glutos''
for the general case.

Glutoses of examples (1)--(4) of section \ref{exa} above are SG-glutoses,
as well as $\cC^{\cD}$ when \cC is an SG-glutos;
the glutos of \'{e}tale schemes is not an SG-glutos.
Any \cU-topos has, canonically, a structure of an
SG-glutos, if one defines open arrows as just those arrows
$u\colon U \arr X$ which locally are mono (here
``locally'' is, of course, with respect to canonical pretopology
of the topos). In fact, the latter example
can be generalized, as shows the following proposition,
easily deduced from ``Giraud theorem'' \ref{gir} and the fact that
\cU-toposes are locally \cU-small (see
p.251 of \cite{SGA4}):

\begin{Prop}
\label{sgglgl}
For any \cU-glutos \cC the presite
${\rm SG}\cC=\cL\cM\cC$
is, in fact, an \cU-glutos.
\end{Prop}

\begin{Rem}
Let \cC be an \cU-presite such that $\cM\cC$
(and, hence, ${\rm SG}\cC$) is an \cU-presite.
Consider the arrow
\begin{equation}
\label{sgarrsg}
\widetilde{{\rm SG}\cC}\arr
{\rm SG}\widetilde{\cC},
\end{equation}
obtained from the universal arrow $Y_{\cC}\colon\cC\arr
\widetilde{\cC}$ by applying the pseudofunctor
$\sim\circ{\rm SG}$ to it and using Prop. \ref{sgglgl}
afterwards. The arrow (\ref{sgarrsg}) is fully faithful
if the pretopology of \cC is subcanonical; the inclusion
arrow (\ref{scharresch}) of sect. \ref{ugps} is the
particular case of the arrow~(\ref{sgarrsg}).
\end{Rem}

The following addition to Theorem \ref{th} states that
the set of SG-presites is stable by the reflection $\sim$:

\begin{Prop}
\label{sggosg}
If \cC is an SG-presite, then the universal glutos for \cC
is an SG-glutos.
\end{Prop}

It is clear from above that any SG-glutos can be obtained
as a universal glutos for some $\cM$-presite \cC
and that the corresponding universal arrow
$Y_{\cC}\colon \cC\arr \widetilde{\cC}$
for an
\cM-\cU-presite \cC can be pulled through the
presite $\cM\widetilde{\cC}$ as
\begin{equation}
\label{unneargl}
\cC\stackrel{Y'_{\cC}}{\arr}
\cM\widetilde{\cC}
\harr\cL\cM\widetilde{\cC},
\end{equation}
where $Y'_{\cC}=\cM Y_{\cC}$.

Call an \cU-presite \cC {\bf nearly \cU-glutos} if
it is naturally equivalent to a presite
$\cM\widetilde{\cC}$ for some \cU-presite \cC.
Meditating over the decomposition (\ref{unneargl}) one
can conclude that the full 2-subcategory of nearly \cU-glutoses
is reflective in that of all \cM-\cU-presites, whereas
the arrow $Y'_{\cC}$ in (\ref{unneargl}) is the unit of the
corresponding adjunction. Besides, the construction
of universal glutos for a nearly \cU-glutos \cC
consists simply in adding of all locally clopen arrows
to the set of clopen arrows.

The next proposition giving an ``internal'' description
of nearly \cU-glutoses is a kind of ``Giraud theorem''
for them.

\begin{Prop}
\label{ngir} Let $(\cC,\tau)$ be an \cM-presite
such that \cC is an \cU-category and
the pretopology $\tau$ is subcanonical. Then
$(\cC,\tau)$ is a nearly \cU-glutos iff the set
$\cO_{\tau}$ of clopen arrows of it satisfies
condition ${\rm(G4_{\cU}+5)}$
as well as the following conditions:

\noindent {\bf(NG1)} For any object $X$ of \cC the set of
clopen subobjects of $X$ is \cU-small;

\noindent {\bf(NG2)} Any family
$\{U_i\stackrel{u_i}{\arr}X\}_{i\in I}$
of clopen arrows has a factorization into a covering
$\{U_i\stackrel{u'_i}{\arr}U.\}_{i\in I}$
followed by a clopen arrow
$U.\stackrel{u.}{\arr}X$ {\rm (in particular,
unions of arbitrary families of clopen
subobjects of $X$ exist (in the lattice of all subobjects
of $X$) and are clopen);}

\noindent {\bf(NG3)} Any epi sink
consisting of clopen arrows is a covering of $\tau$ {\rm (and,
hence, is universal effective epi).}
\end{Prop}

Note that the pretopology of
a nearly \cU-glutos is uniquely recovered from the underlying
clopos structure (just as in the case of \cU-glutoses),
so that we will consider nearly \cU-glutoses either
as presites or as cloposes, depending on circumstances.

\begin{Rem}
\label{mglufun}
The definition of glueing data (see sect. \ref{glu}) with values
in mono's of a category \cC can be essentially
simplified. Namely, define for any set $I$ the category
$\Gamma'I$ as follows. The set of objects of $\Gamma'I$
is the free commutative idempotent monoid
$W(I)/R$ over $I$ (i.e. the set of ``relations''
$R$ consists of two relations: $X^2=X$ and~$XY=YX$).
For any objects $X$ and~$Y$ of $\Gamma'I$ there exists
the only arrow $X \arr Y$ iff there exists
$Z$ such that $YZ=X$. Let $n$ be a natural number. Denote
$\Gamma_nI$ the full
subcategory of $\Gamma'I$ consisting of all monomials
over the variables from $I$ of degree $\leq n$, with
the neutral element of the monoid $W(I)/R$ excluded;
let $\Gamma_+I$ be the union of all $\Gamma_nI$ (it
will be supposed further that the (discrete) category
$\Gamma_1I$ coincides with the set $I$).

Now call a functor $U\colon \Gamma_2I \arr \cC$
with values in mono's of \cC an {\bf \cM-gluon}
or {\bf \cM-glueing data} if there exists its
continuation on $\Gamma_3I$ which respects pullbacks
existing in $\Gamma_3I$ (``cocicle condition''); it then
follows that there exists a continuation $U_+$ of $U$
on the whole $\Gamma_+I$ respecting pullbacks of $\Gamma_+I$
and $U_+$ is unique up to a functor isomorphism.

It is evident enough that one can replace, in the
context of \cM-presites or nearly \cU-glutoses,
open glueing data by ``equivalent'' \cM-glueing data.
\end{Rem}

The following proposition generalizes the realization of
sheaves over topological space $X$ as sheaves of sections
of corresponding fibre bundles over $X$.

\begin{Prop}
\label{nglu}
Let $(\cC,\cO)$ be a nearly \cU-glutos and $\cL\cO$
be the set of all locally clopen arrows of it. Then
for any object $X$ of \cC
the category $\cL\cO/X$ is naturally equivalent
to the SG-topos ${\rm Sh}X$ of sheaves over the complete
Heyting algebra $\cO(X)$ of clopen subobjects of $X$.
\end{Prop}

\noindent{\bf Indications to the proof.} The corresponding
natural equivalence
$J\colon{\rm Sh}X\stackrel{\approx}{\arr}\cL\cO/X$
can be constructed as
follows. Let $s\colon \cO(X)\arr
\cO/X$ be some natural equivalence selecting for any clopen subobject
$u$ of $X$ a clopen arrow $su\colon
U \arr X$ representing this subobject.
Given a sheaf $F\colon \cO(X)\arr {\bf Set}_{\cU}$,
we want to construct a locally clopen
arrow $JF\colon E\arr X$ such that
its sheaf of sections ($u\in \cO(X)\mapsto[su,JF]$) is isomorphic
to $F$. As a first approximation to $JF$ one can take
the coproduct (in $\cL\cO/X$):
\[p=\coprod_{u\in \cO(X)}F(u)\otimes su,\]
where $S\otimes Y$ means the coproduct of the family
$\{Y\}_{i\in S}$ (copower of $Y$). Unfortunately,
$p$ has too many sections as compared to $F$, so that to obtain $JF$
from $p$ one needs to ``glue together'' any two summands of $p$
along the maximal clopen arrow where they are to coincide.

Now we will go from informal considerations above to
the formal constructions.
Define first the ``index set'' $I_F$ as
\[I_F=\coprod_{u\in \cO(X)}F(u);\]
define a partial order relation on $I_F$
such that for any pair $\langle u,x\rangle$,
$\langle v,y\rangle$ ($u,v\in \cO(X)$, $x\in F(u)$,
$y\in F(v)$) of elements of $I_F$ one has that
$\langle u,x\rangle\leq\langle v,y\rangle$
iff $u\leq v$ and $x=\rho^v_uy$, where, of course,
$\rho^v_u\colon F(v)\arr F(u)$
are the corresponding restriction maps of the sheaf $F$.

For any pair $i=\langle u,x\rangle$ and
$j=\langle v,y\rangle$ of elements of $I_F$ there exists
the intersection $i\wedge j=\langle w,z\rangle$,
where $w\leq u\wedge v$ is the biggest element of $\cO(X)$ such that
$\rho^u_wx=\rho^v_wy$ and $z=\rho^u_wx$; note that this
property of $I_F$ essentially depends on the fact that
$F$ is a sheaf and not simply a presheaf.

There exists the only functor
\[\varphi\colon\Gamma_+I_F \arr I_F\]
such that $\varphi$ is the identity map on
$I_F=\Gamma_1I_F\subset\Gamma_+I_F$ and, besides,
for any pair $i$,~$j$ of elements of $I_F$
the identity $\varphi(ij)=i\wedge j$ holds. Recall that the category
$\Gamma_+I$ is defined in Remark \ref{mglufun} above and
that $I_F$ is a category being a partially ordered set.

There is as well an evident forgetful functor
$I_F \arr \cO(X)$ ($\langle u,x\rangle
\mapsto u$), which produces some functor
$N\colon I_F \arr \cL\cO/X$, when being composed
with the chain of functors
\[\cO(X)\stackrel{s}{\arr}\cO/X
\harr\cL\cO/X.\]

Composing now the functor $N$ with the restriction
of the functor $\varphi$ (constructed above) on
the subcategory $\Gamma_2I_F$ of $\Gamma_+I_F$
one obtains some functor
\[U_F\colon\Gamma_2I_F \arr \cL\cO/X.\]
One can verify easily that the functor $U_F$ is
an \cM-gluon, whereas its colimit
$U_F.$ in $\cL\cO/X$ can play the role of the locally
clopen arrow $JF$ corresponding to the sheaf $F$.\eop

\begin{Rem}
In an earlier author's work \cite{Mol} the term ``\cU-glutos''
meant something which is now called ``nearly \cU-glutos'',
whereas \cM-\cU-presites with a subcanonical pretopology
were called there
``\cU-preglutoses''; the main result of \cite{Mol} was, in this terms,
that every \cU-preglutos has a universal completion to an
\cU-glutos, whereas its proof has used generalized
``charts and atlases routine'' (see the next section).
Later on it was observed the presence of SG-toposes ``inside'' glutoses
just via Prop.\ref{nglu}, and the natural question arose
how to generalize both the very notion of glutos
(so that {\it arbitrary\/} toposes can occur in place of
SG-toposes)
and the  theorem of existence of universal glutoses (charts
and atlases method failed to prove Theorem \ref{th} due to
the reasons explained in
Appendix~\ref{Proof}).
{Sect.~\ref{Proof}}
{Sect.~\ref{Proof}}).
\end{Rem}

\section{Charts and Atlases}
\label{charts}
In this section a way of constructing of universal glutoses
(or, rather, of nearly glutoses) by means of charts and atlases
is considered, applicable for \cM-presites with subcanonical pretopology.

Give first some necessary definitions. A continuous functor
$J\colon \cC \arr \cD$ between presites will be said {\bf to
reflect coverings} if for any family
$\{u_i\colon U_i \arr X\}_{i\in I}$
of clopen arrows of \cC the fact that
$\{Ju_i\colon JU_i \arr JX\}_{i\in I}$
is a covering in $\cD$ implies that
$\{u_i\colon U_i \arr X\}_{i\in I}$ is a covering in \cC;
if, besides, \cC is an
\cM-presite with subcanonical pretopology and the underlying functor
of $J$ is faithful then $J$ will be said to {\bf admit atlases}.

An \cM-presite with subcanonical pretopology will be called a
{\bf DG-presite} if it satisfies the factorization condition
(NG2) in Prop.\ref{nglu} above for arbitrary sinks of clopen arrows
(DG above deciphers as ``differential-geometrical'', because presites
of this kind are typical just for differential geometry).

Let \cC be an \cM-\cU-presite with subcanonical
pretopology and $J\colon\cC\arr\cD$ be
an arrow admitting atlases, such that \cD is a nearly
$\cU'$-glutos, where the universe $\cU'$ is any
universe containing \cU as a subset. In constructing
the universal nearly \cU-glutos for \cC one can
use the arrow $J$ considering objects of $\widetilde{\cC}$
as objects of \cD with additional structure.

In fact,
the process of completion of \cC to $\widetilde{\cC}$
using the arrow $J$
can be performed in two steps: first, one completes \cC
with objects which are ``unions of families of clopen subobjects'',
arriving to a universal DG-presite, associated with \cC;
the second step is the completion of the DG-presite so obtained with objects,
which are colimits of clopen \cM-gluons.
Only the second step will be described below,
 because it occurs very frequently in practice.

So let us assume that the arrow
$J\colon \cC\arr \cD$ admitting atlases is given such that
\cC is a DG-\cU-presite, whereas
\cD is a DG-$\cU'$-presite, where the universe
$\cU'$ contains the universe \cU as a subset.

Let $X$ be an object of \cD. An \cU-small  family
$\{U_i\}_{i\in I}$ of objects of \cC together with a
covering $\{JU_i\stackrel{u_i}{\arr}X\}_{i\in I}$
of $X$ will be called an $J${\bf-atlas on} $X$ if
for every $i,j\in I$ the pullback
$JU_i\coprod_XJU_j$ has a representation
\begin{equation}
\label{pbofchart}
\expandafter\ifx\csname xymatrix\endcsname\relax%
\begin{picture}(80,80)
\bx10 0,JU_i;\bx70 0,X;\bx10 60,JU_{ij};\bx70 60,JU_j;
\ar25 5,1 0,40;\ar25 65,1 0,35;
\ar10 55,0 -1,43;\ar70 55,0 -1,43;
\bx0 30,Ju'_i;\bx40 70,Ju'_j;\bx80 30,u_j;\bx40 10,u_i;
\end{picture}
\else
\xymatrix{
JU_{ij}\ar[r]^{Ju'_j}\ar[d]_{Ju'_i}&JU_j\ar[d]^{u_j}\\
JU_i\ar[r]^{u_i}&X
}
\fi
\end{equation}
such that both $u'_i$ and $u'_j$ are clopen arrows of \cC.
Any arrow $u_i$ will be called a {\bf chart} of the corresponding $J$-atlas.

We will identify further a sink
$\{JU_i\stackrel{u_i}{\arr}X\}_{i\in I}$ with a $J$-atlas,
omitting its first component $\{U_i\}_{i\in I}$;
we will write as well ``atlas'' instead of
``$J$-atlas'', when this will not lead to confusion.

\begin{Rem}
The fact that $J$ admits atlases imply that if clopen
arrows $U\stackrel{u}{\arr}V$ and
$U'\stackrel{u'}{\arr}V$ are such that
both $Ju$ and $Ju'$ represent one and the same clopen
subobject of $JV$ then $u$ and $u'$ represent one
and the same subobject of $V$ (i.e. there exists an iso $i$
such that $u'=ui$).

In particular, clopen arrows $u'_i$
in the definition of atlases above (see the pullback (\ref{pbofchart})) are,
essentially, unique, determining,
thus, some clopen \cM-glueing data in \cC such that
$X$ is their ``colimit in \cD''.
\end{Rem}

Given atlases $A$ and $A'$ on $X$ we will say that
$A$ is {\bf compatible with} $A'$ if the union
sink $A\cup A'$ (whose definition is evident)
is an atlas on $X$ as well. One can prove that the relation
between atlases just defined is, in fact, an equivalence
relation; the equivalence class of an atlas $A$
will be denoted further $[A]$.

Let $A=\{JU_i\stackrel{u_i}{\arr}X\}_{i\in I}$
be an atlas on $X$ and
$B=\{JV_k\stackrel{v_k}{\arr}Y\}_{k\in K}$
be an atlas on $Y$. An arrow $f\colon X \arr Y$
will be called $A$-$B${\bf-admissible} if for any
chart $u_i$ of the atlas $A$ and for any chart $v_k$
of the atlas $B$ the pullback of $v_k$ along $fu_i$
has a representation
\begin{equation}
\label{pbadmiss}
\expandafter\ifx\csname xymatrix\endcsname\relax%
\begin{picture}(140,80)
\bx10 0,JU_i;\bx70 0,X;\bx130 0,Y;\bx10 60,JW_{ik};\bx130 60,JV_k;
\ar20 5,1 0,40;\ar80 5,1 0,40;\ar25 65,1 0,95;
\ar10 55,0 -1,43;\ar130 55,0 -1,43;
\bx-3 30,Jw_{ik};\bx70 70,Jf_{ik};
\bx40 10,u_i;\bx100 10,f;\bx140 30,v_k;
\end{picture}
\else
\xymatrix{
JW_{ik}\ar[rr]^{Jf_{ik}}\ar[d]_{Jw_{ik}}&&JV_k\ar[d]^{v_k}\\
JU_i\ar[r]^{u_i}&X\ar[r]^f&Y
}
\fi
\end{equation}
such that $w_{ik}$ is a clopen arrow of \cC.

\begin{Prop}
\label{correct}
If an arrow $f\colon X\arr Y$ of \cD is
$A$-$B$-admissible for some atlases $A$ and $B$, then
$f$ is $A'$-$B'$-admissible for any atlases $A'\in[A]$
and $B'\in[B]$; if, besides, an arrow
$g\colon Y \arr Z$ is $B$-$C$-admissible,
then the composition arrow $gf$ is $A$-$C$-admissible.
\end{Prop}

The latter proposition justifies correctness of the following
definitions and constructions. First, call
the arrow $f$ above $[A]$-$[B]$-{\bf admissible} if it is
$A$-$B$-admissible. Now one can define the category
$\cC_J$ as follows. Objects of $\cC_J$ are
all pairs $\langle X,[A]\rangle$ consisting of an
object $X$ of \cD and an equivalence class $[A]$
of atlases on it. Arrows of $\cC_J$ are all triples
$\langle\langle X,[A]\rangle,f,\langle Y,[B]\rangle\rangle$
such that the arrow $f\colon X \arr Y$ is
$[A]$-$[B]$-admissible (and we will write simply $f$ instead of
the whole triple in situations not leading to confusion).

There are evident functors $J_{\cC}\colon\cC
\arr \cC_J$ ($X\mapsto\langle JX,[\{Id_X\}]\rangle$) and
$J'\colon \cC_J \arr \cD$ (forget
atlases). There is the natural pretopology on
$\cC_J$ making both $J_{\cC}$ and $J'$
continuous. This pretopology is defined as follows.
Declare a monic arrow $f\colon X \arr Y$ between
objects $\langle X,[A]\rangle$ and $\langle Y,[B]\rangle$
of $\cC_J$
$J${\bf-clopen} if all arrows $f_{ik}$ in the diagram (\ref{pbadmiss})
above are clopen arrows of \cC (it follows then from the condition (NG2)
that $f$ is a clopen arrow of \cD).
Let $\tau$ consists
of all sinks $S$ in $\cC_J$ such that any arrow of
$S$ is $J$-clopen and $J'S$ is a covering in \cD.

One can prove that $\tau$ is really a pretopology
on $\cC_J$ and the functors $J_{\cC}$ and~$J'$
become continuous if one equips the category $\cC_J$
with the pretopology $\tau$. The notations $\cC_J$,
$J_{\cC}$ and~$J'$ will be reserved as well to denote
the corresponding presite and morphisms of presites.
Note that the equality $J=J'J_{\cC}$ holds.

Now, at last, one can formulate the theorem giving
a construction of universal nearly \cU-glutos
by means of charts and atlases.

\begin{Th}
\label{thcharts}
Let \cC be a {\rm DG}-\cU-presite, \cD be
a {\rm DG}-$\cU'$-presite for $\cU\subset \cU'$.
Let the arrow $J\colon \cC\arr \cD$ admits
atlases. Then the presite $\cC_J$ constructed above
is a {\rm DG}-\cU-presite. If, moreover, \cD
is a nearly $\cU'$-glutos then $\cC_J$ is a universal
nearly \cU-glutos for \cC, whereas the arrow
$J_{\cC}\colon \cC\arr \cC_J$ is
a corresponding universal arrow.
\end{Th}

Applying this theorem to standard constructions of
differential geometry (manifolds, vector bundles,
principal $G$-bundles, etc.,) one can check that all this
constructions are just particular cases of universal
(nearly) \cU-glutos construction. But to check that
certain functors of algebraic geometry like Spec above
fall as well into this scheme, one needs another tools.
The theorem below gives sufficient criteria for an
arrow between presites to be universal.

Before formulating this theorem one needs to
introduce one more definition. An arrow
$F\colon \cC\arr \cD$
will be said to {\bf locally reflect clopens}
if for any arrow $u\colon U \arr X$ of \cC
the fact that $Fu$ is clopen implies that $u$ is locally clopen.

%
\begin{Th}
\label{th3}
Let \cC be an \cM-presite with a subcanonical pretopology,
\cD be a nearly \cU-glutos and $Y\colon \cC\arr \cD$
be a continuous functor. Then the following statements are equivalent:

\noindent {\rm (a)} $Y$ is a universal arrow for \cC;

\noindent {\rm (b)} $Y$ is fully faithful, reflects coverings,
locally reflects clopens and, besides,
for every object $D$ of \cD there exists an \cU-small
covering $\{u_i\colon YU_i \arr D\}_{i\in I}$ of
$D$ by ``objects of \cC''.
\end{Th}

Now the universality of the arrow ${\rm Spec}$ can be
established just with the help of Theorem \ref{th3}.
This theorem can be applied as well to obtain necessary
and sufficient conditions for a given \cU-valued
functor on the category {\bf Ring} to be representable
by a Grothendieck scheme; these conditions can be formulated
in terms of Zariski pretopology on the category ${\bf Ring}^{op}$.
(cf. the  existence problem of Grothendieck as formulated in \cite{Mum}).

\appendix
\def\la{\langle}
\def\ra{\rangle}
\ifx\Sh\undefined\rmname{Sh}\fi
\oplim{colim}
\mathdef\Tsub$T_{sub}$
\mathdef\Csub$\cC_{sub}$
\mathdef\Osub$\cO_{sub}$
\mathdef\tsub$\tau_{sub}$
\mathdef\vd$\vdash$
\mathdef\vD$\vDash$
\mathdef\Shu$\Sh_{\cU'}\cC$
\vers
{\section{The idea of the proof of Th.\protect\ref{th}}
\label{Proof}}
{\section{The Proof of Th.\protect\ref{th}}
\label{Proof}
\subsection{The scetch of the proof of Th.\protect\ref{th}}
\label{proof}}
{\section{The Proof of Th.\protect\ref{th}}
\label{Proof}
\subsection{The scetch of the proof of Th.\protect\ref{th}}
\label{proof}}
Let a set $\cU'$ be
a universe such that $\cU\subset \cU'$ and \cC
is  $\cU'$-small (recall that we are
living, due to Sect. \ref{sets}, in ``Grothendieck's paradise'' restricted from above by the universal class of all sets).
Let
${\rm Sh}_{\cU'}\cC$ be the topos of $\cU'$-valued
sheaves on \cC, considered as a presite via canonical
pre\-to\-po\-logy~$\tau$.

In constructing the universal arrow $Y_{\cC}\colon \cC
\arr \widetilde{\cC}$ the Yoneda functor
$Y\colon \cC\arr
{\rm Sh}_{\cU'}\cC$ can be used, whereas
$\cU'$-valued sheaves can be considered as building
blocks in the process of construction of $\widetilde{\cC}$.

In more details, let $T$ be a (non-elementary) theory
whose axioms are axioms of elementary theory of categories
together with axioms (PT1)--(PT4) of presites and conditions
(G2), ${\rm(G4_{\cU})}$, (G5) and ${\rm(G6_{\cU})}$
(as well as (G5P) in case of weak form of (G5))
imposed on the set of clopen arrows (these conditions are the same as
in ``Giraud theorem'' \ref{gir}, {\it excepting\/}
the \cU-smallness condition ${\rm(G7_{\cU})}$). The presite
${\rm Sh}_{\cU'}\cC$ is easily seen to be a model
of the theory $T$ and one can prove that submodels of $T$
form a complete lattice with respect to the inclusion
functors. The latter lattice is, essentially,
a closure system on the set
\[ X={\rm Mor}({\rm Sh}_{\cU'}\cC)\coprod\tau\ .\]

One can prove that the $T$-closure of the image $Y\cC$
of \cC by Yoneda functor in ${\rm Sh}_{\cU'}\cC$
is not only model of $T$ but satisfies the condition
${\rm(G7_{\cU})}$ as well. In other words, it is an \cU-glutos
and one can show further that it is the universal glutos
$\widetilde{\cC}$.

In proving this it is useful to ``translate'' axioms of
the theory $T$ into the set of {\it inference rules\/}
on the set $X$
\vers
{(in the sense of \cite{Acz}),}
{(see Appendix~\ref{clopf} for definitions of formal systems,
inference rules, etc. on complete lattices;
our definitions
are a bit more general than those
given in~\cite{Acz})),}
{(see Appendix~\ref{clopf} for definitions of formal systems,
inference rules, etc. on complete lattices;
our definitions
are a bit more general than those
given in~\cite{Acz})),}
whereas
arrows and coverings in $Y\cC$ to consider as {\it axioms\/}
of the corresponding (infinitary) formal system (denoted
further $FS(T)$). Then the $T$-closure of $Y\cC$
in ${\rm Sh}_{\cU'}\cC$ turns out to be, essentially,
the set of {\it theorems\/} of the formal system $FS(T)$.

It is convenient (as well as more informative) to separate
the subtheory $T_{sub}$
of ``presites with subcanonical pretopology'' in $T$;
considering first the $T_{sub}$-closure of $Y \cC$
one can prove that the full sub-2-category of subcanonical
\cU-presites is reflective in \mbox{{\bf Psite}${}_{\cU}$}. This reduces the proof of Theorem \ref{th} to
the particular case of \cU-presites \cC
with subcanonical pretopology.

In proving that both $T_{sub}$-closure $\cC_{sub}$
and $T$-closure $\widetilde{\cC}$ of
$YC$ in ${\rm Sh}_{\cU'}\cC$ have \cU-small
local sets of topological generators (see condition
(${\rm LTG}_{\cU}$) of sect. \ref{ps}) the following Lemma,
easily deduced from Lemme 3.1 on p.231 of \cite{SGA4},
\vers
{is crucial:}%
{is crucial (in what follows we will often write
$\wt X$, etc. instead of $YX$, etc.):}%
{is crucial (in what follows we will often write
$\wt X$, etc. instead of $YX$, etc.):}%

\begin{Lem}
\label{Lem0}
Let \cC be an \cU-presite. Then the Yoneda map
$Y\colon\cC\arr
{\rm Sh}_{\cU'}\cC$ has the following property:
\vers
{for any objects $X$ and $X'$ of \cC and any arrow
$f\colon YX \arr YX'$ there exists a covering
$\{V_i\stackrel{v_i}{\arr}X\}_{i\in I}$ in \cC
such that the set of objects $\{V_i:i\in I\}$ is a subset
of the set $G_X$ of
topological
generators over $X$ and for any
$i\in I$ there exists an arrow $v'_i\colon V_i \arr X'$
such that the identity $f\comp Yv_i=Yv'_i$
holds.}%
{for any objects $X$ and $X'$ of \cC and any arrow
$f\colon\wt X\arr\wt{X'}$ there exists a covering
$\{v\colon V_v\arr X\}_{v\in V}$ of $X$ in \cC
such that the set of objects $\{V_v:v\in V\}$ is a subset
of the set $G_X$ of
topological
generators over $X$ and for any
$v\in V$ there exists an arrow $f_v\colon V_v \arr X'$
such that the identity $f\comp\wt v=\wt{f_v}$
holds:}%
{for any objects $X$ and $X'$ of \cC and any arrow
$f\colon\wt X\arr\wt{X'}$ there exists a covering
$\{v\colon V_v\arr X\}_{v\in V}$ of $X$ in \cC
such that the set of objects $\{V_v:v\in V\}$ is a subset
of the set $G_X$ of
topological
generators over $X$ and for any
$v\in V$ there exists an arrow $f_v\colon V_v \arr X'$
such that the identity $f\comp\wt v=\wt{f_v}$
holds:}%
\ifnum\version>1
\begin{equation}
\label{X-cover}
\xymatrix@1{
{\wt{V_v}}\ar[r]^{\wt v}\ar@/_1pc/[rr]_{\wt{f_v}}&{\wt{X}}\ar[r]^f&{\wt{X'}}
}
\end{equation}
The family
$\{v\colon V_v\rightarrow X,f_v\colon V_v\rightarrow X'\}_{v\in V}$
will be called a {\bf $G_X$-covering of the arrow $f$}.
\fi
\end{Lem}

It is just applications of this lemma in transfinite induction on the
length of {\it proofs\/} (in formal systems
$FS(T_{sub})$ and $FS(T)$) which permits one to prove
that both $\cC_{sub}$ and $\widetilde{\cC}$ are
\cU-presites. Moreover, one can prove that any
object, arrow and covering of the $T_{sub}$-closure
$\cC_{sub}$ has a {\it finite\/} proof, which permits one
to describe the presite $\cC_{sub}$ explicitly.

Now the universality properties of the corresponding
arrows $Y'_{\cC}\colon \cC\arr
\cC_{sub}$ and $Y_{\cC}\colon \cC\arr
\widetilde{\cC}$ follow from that of ``sheafified
Yoneda functors'' if one applies (transfinite) induction on the length of proofs: given a continuous functor
$F\colon\cC\arr \cD$ into a subcanonical
\cU-presite, resp. into an \cU-glutos one has that
any element $Z$ (arrow or covering) of $\cC_{sub}$, resp.
of $\widetilde{\cC}$ having a proof $P$, where
a family of axioms $\{A_i\}_{i\in I}$ from \cC were used,
goes by the functor
\[{\rm Sh}_{\cU'}(F)\colon {\rm Sh}_{\cU'}(\cC)
\arr {\rm Sh}_{\cU'}(\cD)\]
into an element $Z'$ which has ``the same'' proof
in ${\rm Sh}_{\cU}(\cD)$ as $Z$ has in
${\rm Sh}_{\cU}(\cC)$ with only the family
$\{A_i\}_{i\in I}$
of axioms replaced by the family $\{FA_i\}_{i\in I}$.
This implies that $Z'$ belongs to the closure of $\cD$
(naturally equivalent to \cD, because \cD is a model
of $T_{sub}$, resp. of $T$),
i.e.
the restriction of the functor
${\rm Sh}_{\cU}(F)$ on $\cC_{sub}$, resp. on
$\widetilde{\cC}$ can be pulled through $\cD$.

%
It turns out, that if \cC is an SG-presite with subcanonical
pretopology, then any theorem of the formal
system $FS(T)$ has a proof of a fixed finite length.
In this case one can use as well another continuous functors
$F\colon\cC\arr\cD$ in place of Yoneda functor in constructing of
$\widetilde{\cC}$
(namely, functors admitting atlases defined in sect. \ref{charts} above).
\ifnum\version>1
\subsection{Construction of the Universal Subcanonical Presite}
\label{proof-sub}
Construct now the universal arrow
\begin{equation}
\label{Y-sub}
Y_{sub}\colon\cC\arr\cC_{sub}
\end{equation}
from an \cU-presite \cC to a presite with subcanonical pretopology.

We are to prove, first of all, that for any presite \cC with subcanonical
pretopology (i.~e. a model of the theory $T_{sub}$) its subpresites
with subcanonical pretopology (i.e. submodels of the theory $T_{sub}$)
form a closure system in the complete lattice of all subpresites of \cC
(see Appendix~\ref{clopf}, especially Example~\ref{fE4}).
I.~e., we are to prove that the intersection of any family of subpresites
with subcanonical pretopology (which {\it is\/} a subpresite due to
Example~\ref{fE4} of Appendix~\ref{clopf}) has subcanonical pretopology
itself.

Let $\{\cC_i\}\iI$ be a family of subpresites of \cC with subcanonical
pretopology. Let a family $\{\sar Uu,u,X;\}\uU$ be a covering in
$\bigcap_i\cC_i$.
For any $i\in I$ one has
\[
\expandafter\ifx\csname xymatrix\endcsname\relax%
\divide\dgARROWLENGTH by4
X=
\colim_{\cC_i}
\left\{
\begin{diagram}
\\
\node[2]{U_u}\\
\node{U_{u}\prod\nolimits_XU_{u'}}\arrow{ne}\arrow{se}\\
\node[2]{U_{u'}}\\
\end{diagram}
\right\}
=
\colim_\cC
\left\{
\begin{diagram}
\\
\node[2]{U_u}\\
\node{U_u\prod\nolimits_XU_{u'}}\arrow{ne}\arrow{se}\\
\node[2]{U_{u'}}\\
\end{diagram}
\right\},
\else
X=
\colim_{\cC_i}
\left\{
\vcenter{
\xymatrix@1@=3pt{
&&&U_u\\
U_{u}\prod\nolimits_XU_{u'}
\ar[urrr]\ar[drrr]\\
&&&U_{u'}
}}
\right\}
=
\colim_\cC
\left\{
\vcenter{
\xymatrix@1@=3pt{
&&&U_u\\
U_{u}\prod\nolimits_XU_{u'}
\ar[urrr]\ar[drrr]\\
&&&U_{u'}
}}
\right\},
\fi
\]
because of the continuity of inclusion of subpresite functors.

We are to prove that
\[
X=
\colim_{\bigcap_i\cC_i}
\expandafter\ifx\csname xymatrix\endcsname\relax%
\left\{
\divide\dgARROWLENGTH by4
{\begin{diagram}
\\
\node[2]{U_u}\\
\node{U_{u}\prod\nolimits_XU_{u'}}\arrow{ne}\arrow{se}\\
\node[2]{U_{u'}}\\
\end{diagram}}
\right\}.
\else
\left\{
\vcenter{
\xymatrix@1@=3pt{
&&&U_u\\
U_{u}\prod\nolimits_XU_{u'}
\ar[urrr]\ar[drrr]\\
&&&U_{u'}
}}
\right\}.
\fi
\]

But, given a family $\{\sar U_u,f_u,Y;\}\uU$ in $\bigcap_i\cC_i$
such that $f_u|_{U_u\prod_XU_{u'}}=f_{u'}|_{U_u\prod_XU_{u'}}$,
for any $i\in I$ there exists in $\cC_i$ the only $f\colon X\arr Y$
glued out of $f_i$; this $f$ is as well the only in \cC, because
$\cC_i\subset\cC$ is continuous. Hence, this $f$ belongs to
$\bigcap_i\cC_i$, i.~e. the family
$\{\sar Uu,u,X;\}\uU$
is effective epi in
$\bigcap_i\cC_i$, which implies that it is universal effective epi family.
This proves that {\bf sub-$T_{sub}$-models of a subcanonical presite
form a closure system.}

What is the minimal set of ``schemes of inference rules'' for the theory
\Tsub? Or, rather, what inference rules are to be added to the
inference rules (C)--(PT4$''$) of Example~\ref{fE4} of Appendix~\ref{clopf}
to get inference rules determining \Tsub? Our guess is that it is sufficient
to add the following scheme:
\begin{list}{}{\topsep1mm\itemsep1mm}
\item[\bf(S)] Let
$\la X,U\ra$
be a covering of $\tau$ and $f\colon X\arr Y$.
Then:
\begin{equation}
\label{S-rule}
\la X,U\ra, \{f|_{\partial_0u}:=f\comp u\}\uU\vd f
\end{equation}
\end{list}

Really, it is evident, that any submodel
$\la\cC',\tau'\ra\subset\la\cC,\tau\ra$ of \Tsub-model must be closed
with respect to the inference rules~(S);
one needs to prove only that any subpresite
$\la\cC',\tau'\ra$ of the presite $\la\cC,\tau\ra$ which is closed
with respect to the inference rules~(S) is a \Tsub-submodel of
$\la\cC,\tau\ra$. In other words, one needs to prove that the pretopology
$\tau'$ is subcanonical, i.e. that any covering
$\{u\colon U_u\arr X\}\uU$ of $X$ in $\tau'$ is effective epi in $\cC'$,
i.e. that
\[
X=
\colim_{\cC'}
\left\{
\vcenter{
\xymatrix@1@=3pt{
&&&U_u\\
U_{u}\prod\nolimits_XU_{u'}
\ar[urrr]\ar[drrr]\\
&&&U_{u'}
}}
\right\}.
\]

Really, let $Y\in\cC'_0$ and for any $u\in U$ there is given an arrow
$f_u\colon U_u\arr Y$ of $\cC'$ such that
$f_u|_{U_{uu'}}=f_{u'}|_{U_{uu'}}$ (here $U_{uu'}:=U_u\prod_X U_{u'}$).
The fact that $\{u\}\uU$ is a covering in $\tau$
(because
$\la\cC',\tau\ra\subset\la\cC,\tau\ra$
is continuous) implies that there exists the only $f$ in $\cC$ such that
for any $u\in U$ the identity $f|_{U_u}=f_u$, because $\tau$ is subcanonical.
Now $f$ belongs to $\cC'$ by inference rule (S) and this is the only such
$f$. Hence, $\tau'$ is subcanonical, which proves our guess that
the inference rules (C)--(PT4$''$)
of Example~\ref{fE4} of Appendix~\ref{clopf} together with inference rules
(S) determine just \Tsub-submodels of $\la\cC,\tau\ra$.

Let us try to describe now the \Tsub-closure
$\la\cC_{sub},\tau_{sub}\ra$ in \Shu of $\la Y\cC,Y\tau\ra$, where
$Y\colon\cC\arr\Shu$ is Yoneda functor (combined with associated sheaf
functor).


Call an open arrow $u\colon U\arr X$ of \tsub {\bf plain} if there exists
an arrow $f\colon X\arr\wt Y$ in \Csub and an open $v\colon V\arr Y$ in
\cC such that $u$ is a pullback of $\wt v\colon\wt V\arr\wt Y$ along $f$:
\[
\xymatrix{
U\ar[r]\ar[d]_u&{\wt V}\ar[d]^{\wt v}\\
X\ar[r]^f&{\wt Y}
}
\]
Call $u$ {\bf semiplain} if it has a decomposition $u=u_1\cdots u_n$, where
every $u_i$ is plain.

An object $X$ of \Csub will be said to {\bf have a tail} if there exists
a semiplain arrow $f\colon X\arr\wt Y$ for some object $Y$ of \cC; the
arrow $f$ itself will be called a {\bf tail} of $X$.
\begin{Prop}
\label{Proof1}
One has:
\begin{list}{}{\itemindent-\parindent\topsep1mm\itemsep1mm}
\item[\bf a)] Every object $X$ of \Csub has a tail;
\item[\bf b)] Every open arrow $u\colon U\arr X$ of \Csub is semiplain;
\item[\bf c)] \Csub is a full subcategory of \Shu;
\item[\bf d)] If $f\colon X\arr\wt Y$ is a tail of an object $X$ of \Csub,
and $G_Y$ is any set of topological generators in $\cO/Y$,
then the set $G_X=\wt{G_Y}$ is a set of topological generators
in $\Osub/X$.
{\rm (in particular, $G_X$
is \cU-small)}.

Moreover, for any open arrow $u\colon U\arr X$ in \Csub
a covering
$\{v\colon\wt V_u\arr U\}_{v\in V}$
of $U$ by objects of $\wt{G_Y}$ can be chosen in such a way
that for any $v\in V$ one has $fuv=\wt{f_v}$ for some open
arrow $f_v\colon V_v\arr Y$ of \cC
{\rm(such a covering will be called a {\bf $G_Y$-covering}
of an open arrow $u$).}
\end{list}
\end{Prop}
\noindent{\bf Proof.}
To prove a) and b) it is sufficient to see that the set of objects
with a tail, resp. of semiplain open arrows is closed with respect to
the inference rules (C)--(S). But this is evident, because the only
inference rules producing new {\it objects\/} or {\it open arrows\/}
from old ones are (C) and the rules~(\ref{f60a}) of (PT2$'$).

To prove c) and d) we need some more definitions. Define the {\bf height}
of an open arrow $u\colon U\arr X$ of \Csub as the smallest number
of arrows in the representation
$u=(U\rightarrow X_1\rightarrow\cdots\rightarrow X_n=X)$
of $u$ as a composition of plain arrows%
.
The {\bf height} of
an object $X$ of \Csub define as the smallest number $n\in\N$ such that
$X$ has a tail $u\colon X\arr\wt Y$ of height n, or zero if $X=\wt{X'}$
for some object $X'$ of \cC. At last, the {\bf biheight} of an
{\bf open arrow} in \Csub
$u\colon U\arr X$ define as the pair $(n,m)\in\N^+\times\N$,
where $n$ is the height of $u$,
and $m$ is the height of $X$; the {\bf biheight} of an arbitrary {\bf arrow}
$f\colon X\arr Y$ in \Csub define as the pair $(n,m)\in\N\times\N$,
where $n$ is the height of $X$
and $m$ is the height of $Y$.

\ifnum\version=3
One sees immediately that the definitions above determine correctly
some \N-grading
\[
(\Csub)_0=\coprod_{i\in\N}(\Csub)_0^i
\]
on the set of objects of \Csub, as well as bigradings
\[
\Osub=\coprod_{(i,j)\in\N^+\times\N}\Osub^{ij}
\]
and
\[
(\Csub)_1=\coprod_{i,j\in\N}(\Csub)_1^{ij}
\]
on the set of open arrows, resp. on the set of all arrows of \Csub.
\fi

To prove c) and d) we will use induction on biheight. First of all,
the following corollaries of Lemma~\ref{Lem0} prove c) for arrows
of biheight $(0,0)$ as well as d) for open arrows of biheight $(1,0)$:
\begin{Cor}
\label{cor1}
For any objects $X,Y$ of \cC any arrow $f\colon\wt X\arr\wt Y$ belongs to
$\Csub$.
\end{Cor}
\noindent{\bf Proof.}
Let a family
$\{v\colon V_v\rightarrow X,f_v\colon V_v\rightarrow Y\}_{v\in V}$
be a $G_X$ covering of the arrow $f$ (see Lemma~\ref{Lem0}). Then
the family
$\{\wt{v}\colon\wt{V_v}\arr\wt X\}_{v\in V}$ is to be a covering of $\wt X$
in $\Csub$ by continuity of the shifified Yoneda functor. The identities
$f|_{\wt{V_v}}:=f\wt v=\wt{f_v}$ imply that $f|_{\wt{V_v}}\in\Csub$,
which implies, by inference rules (S) that $f\in\Csub$.
\eop
\begin{Cor}
\label{cor2}
For any $f\colon\wt X\arr\wt Y$ and any open arrow $u\colon U\arr Y$
of \cC, there exists a covering
$\{v\colon\wt V_v\arr f^*\wt U\}_{v\in V}$ in \Csub such that the set
$\{\,V_v:v\in V\,\}$
is a subset of the set $G_X$ of topological generators of $X$.
\end{Cor}
\noindent{\bf Proof.}
First of all, let a family
$\{v_i\colon V_i\arr X, f_i\colon V_i\arr Y\}\iI$
be a $G_X$-covering of $f$
(see Lemma~\ref{Lem0}).
Consider the following picture of pullbacks in \Shu:
\[
\xymatrix{
{\wt{V_{ij}}}\ar@{-->}[r]&{\wt{V'_i}}\ar[r]\ar[d]_{\wt{f_i^*u}}\ar@{-->}[dr]
\ar@{-->}@/^1pc/[rr]^{\wt{u^*f_i}}
&f^*\wt U\ar[r]\ar[d]&{\wt U}\ar[d]^{\wt u}\\
&{\wt{V_i}}\ar[r]^{\wt{v_i}}\ar@{-->}@/_1pc/[rr]_{\wt{f_i}}
&{\wt X}\ar[r]^f&{\wt Y}
}
\]
Taking into account that the ``composite '' pullback in the diagram
above is the image of the pullback in \cC of $u$ along $f_i$
(because shifified Yoneda functor respects finite limits), one concludes
that diagonal dashed arrow $\wt{V'_i}\arr\wt{X}$ above is to be open, and
applying Lemma~\ref{Lem0} once again we will find, for any $i\in I$,
a covering $\{V_{ij}\arr V'_i\}_{j\in J_i}$ such that any $V_{ij}$ belongs to
the set $G_X$ of topological generators of $\cO/X$.

Now the family
$\{\wt{V_{ij}}\arr\wt{V'_i}\arr f^*\wt U\}$
is a covering of $f^*\wt U$ in \Csub by inference rules (PT3$'$) and
gives the desired covering of $f^*\wt U$ by objects of $G_X$.
\eop

Prove now d) for open arrows of biheight $(n+1,0)$ as well as c) for
arrows of biheight
$(n+1,0)$
supposing the corresponding statements are valid for (open) arrows
of biheight
$(n,0)$.

Let an open arrow $u'$ in \Csub be a composition
$u'=\sar f^*\wt U,f^*\wt u,{};\sar Z,z,\wt X;$
of an open arrow $z\colon Z\arr X$ of biheight $(n,0)$ and
and a pullback of an open arrow $\wt u\colon\wt U\arr\wt Y$
along an arrow $f\colon Z\arr\wt Y$ ($f$ belongs to \Csub due to
inductive assumption).

There exists, by inductive assumption, an open family
$\{v'_i\colon V_i\arr X\}\iI$ such that the family $\{\wt{v'_i}\}\iI$
factors through $z$ as $\wt{v'_i}=zv_i$ and the family
$\{v_i\}\iI$ is a covering of $Z$ in \Csub.

Now choose, for any $i\in I$ some  $G_Y$-covering
$\{\sar{V_{ij}},v_{ij},{V_i};,\sar{V_{ij}},{f_{ij}},Y;\}_{j\in J_i}$
of the arrow $fv_i$. It exists by Lemma~\ref{Lem0}.

Taking a pullback of $u$ along $f_{ij}$
and applying shifified Yoneda functor $\wt{*}$ afterwards,
one obtains the following commutative picture
\[
\xymatrix{
{\wt{V_{ijk}}}\ar@{-->}[r]&{\wt{f_{ij}{}^*U}}
\ar[rr]\ar[d]_{\wt{f_{ij}{}^*u}}
\ar@{-->}@/^1pc/@<1ex>[rrr]^{\wt{u^*f_{ij}}}
&&f^*\wt U\ar[r]\ar[d]^{f^*\wt u}&{\wt U}\ar[d]^{\wt u}\\
&{\wt{V_{ij}}}\ar[r]^{\wt{v_{ij}}}
\ar@{-->}@/_1pc/@<-1ex>[rrr]_<(.3){\wt{f_{ij}}}
&{\wt{V_i}}\ar[r]^{v_i}\ar@{-->}[dr]_{\wt{v'_i}}
&Z\ar[r]^f\ar[d]^z&{\wt Y}\\
&&&{\wt X}
}
\]

Noting that the arrow
$(f_{ij}{^*}U\arr X)=v'_i\comp v_{ij}\comp f_{ij}{}^*u$ is open in \cC,
we obtain a covering
$\{V_{ijk}\arr{f_{ij}{}^*U}\}$ by elements of $G_X$
(because \cC satisfies (LTG${}_\cU$), being an \cU-presite).
Applying Yoneda gives a covering
$\{\wt{V_{ijk}}\arr\wt{f_{ij}{}^*U}\}$ by elements of $\wt{G_X}$.
Composing the latter covering with the covering
(by inference rules (PT2$'$))
$\{(f^*\wt u)^*(v_i\wt{v_{ij}})\colon\wt{f_{ij}{}^*U}\arr f^*\wt U\}$
and applying inference rules (PT3$'$)
we prove d) for all open arrows over objects of the kind $\wt X$:
the family $\{\wt{V_{ijk}}\arr f^*\wt U\}$ is the desired covering.
It remains only to prove c) for any arrow
$g\colon f^*\wt U\arr\wt{Y'}$, but, composing $g$ with the arrow
$\wt{V_{ijk}}\arr f^*\wt U$ and applying Corollary~\ref{cor1}
we have that for any $i,j,k$ the restriction
$g|_{\wt{V_{ijk}}}$ belongs to \Csub, hence $g$ itself belong to \Csub
by inference rules (S).

Now, in fact, d) is proved for any opens: take an open $u\colon U\arr X$
and augment $X$ by a tail $\sar U,u,X;\sar{},x,\wt Y;$ applying just proved
to an open $xu$.

The only thing which remains to prove is c) for arrows of arbitrary
biheight.

Let $X,Y$ be objects of \Csub and $f\colon X\arr Y$ be some arrow in \Shu.
Consider a commutative diagram of pullbacks in \Shu
\[
\xymatrix{
W_{ij}\ar[r]^{w_{ij}}\ar[d]&f^*\wt{U_i}\ar[r]^{f_i}\ar[d]&{\wt{U_i}}
\ar[d]_{u_i}
\\
{\wt{V_j}}\ar[r]^{v_j}
&X\ar[r]^f\ar[d]^x&Y\ar[d]_y\\
&{\wt{X'}}&{\wt{Y'}}
}
\]
where $x$ and $y$ are tails of $X$, resp. of $Y$,
$\{u_i\colon\wt{U_i}\arr Y\}\iI$
is a $G_{Y'}$-covering of Y,
$\{v_j\colon\wt{V_j}\arr X\}_{j\in J}$
is a $G_{X'}$-covering of X
and the pullback along $f$ is taken in \Shu; the second square is
a pullback in \Shu. We are to prove that the arrow $f$ belongs to \Csub.

Suppose first that for any $i,j$ the object $W_{ij}$ belongs to \Csub
and the arrow $W_{ij}\arr\wt{V_j}$ is open there. Then, applying c)
for arrows of biheight $(m,0)$ one sees that the arrow
$f_iw_{ij}$ will belong to \Csub. And the commutativity of the above
diagram implies that the restriction $f|_{W_{ij}}=u_if_iw_{ij}$
on the (open!) arrow $W_{ij}\arr X$ belongs to \Csub, hence, $f$ itself
belongs to \Csub by inference rules (S).

The following Lemma shows that the arrows $W_{ij}\arr\wt{V_j}$ always
belong to \Csub and are open there.

\begin{Lem}
\label{Lem1}
Consider a pullback diagram
\[
\xymatrix{
g^*U\ar[r]^h\ar[d]_{g^*u}&U\ar[d]^u\\
{\wt V}\ar[r]^g&Y
}
\]
in \Shu, where $u$ is an open arrow in \Csub, and $g\colon\wt V\arr Y$ is
an arbitrary arrow of \Shu. Then the arrow
$g^*u\colon g^*U\arr\wt V$ belongs to \Csub and is open there.
\end{Lem}
\noindent{\bf Proof.}
Use induction by the height of $u$. If height of $u$ is $1$,
i.e. $u=g'^*\wt{u'}$ for some open arrow $u'\colon U'\arr Y'$ of \cC
and some arrow $g'\colon Y\arr\wt{Y'}$ of \Csub
then,
considering the pullback diagram
\[
\xymatrix{
g^*U\ar[r]^h\ar[d]_{g^*u}&U\ar[d]^u\ar[r]&{\wt{U'}}\ar[d]^{\wt{u'}}\\
{\wt V}\ar[r]^g&Y\ar[r]^{g'}&{\wt{Y'}}
}
\]
one concludes, that $g'g$ belongs to \Csub by Cor.~\ref{cor1}, hence,
$g^*u$ belongs to \Csub by inference rules (PT2$'$) and is open in \Csub.

Supposing that Lemma~\ref{Lem1} is true for open arrows of height $n$,
prove it for open arrow $u$ of height $n+1$. Let
$u=(\sar U,u',Z;\sar{},u'',Y;)$, where $u''$ is of height $n$ and $u'$
is of height $1$.

Consider the diagram of pullbacks:
\[
\xymatrix{
g'^*U\ar[r]\ar[d]_{g'^*u'}&U\ar[r]\ar[d]^{u'}&{\wt{U'}}\ar[d]^{\wt{u'''}}\\
g^*Z\ar[r]^{g'}\ar[d]_{g^*u''}&Z\ar[r]^{g''}\ar[d]^{u''}&{\wt{Y'}}\\
{\wt V}\ar[r]^g&Y
}
\]
By inductive assumption, $g^*u''$ is open arrow of \Csub. Hence,
$g^*Z$ belongs to \Csub, which implies that $g''g'\in\Csub$ by c) for
arrows of biheight $(m,0)$ (proved above). Hence, $g'^*u'\in\Csub$
and is open there by (PT2$'$). Now $g^*u=g^*u''g'^*u'$ belongs to \Csub
by inference rules (C) and is, evidently, open there.
\eop\ %

Proposition~\ref{Proof1} is proved.
\eop

Proposition~\ref{Proof1} and its proof shows, in particular, that
any object, arrow or open arrow of \Csub has a finite proof from
axioms $Y\cC$ and $Y\cO$ (one only needs to note that for
any object or open arrow of \Csub there exists a proof using
coverings of a very specific type: $\la X,\{X,u\colon U\arr X\}\ra$,
containing sinks of only two arrows, one of which is the identity arrow).

In fact, to construct the clopos $\la\Csub,\cO_{sub}\ra$ alone,
we can {\it exclude} coverings from inference rules (C)--(PT4$''$)
replacing simultaneously the inference rule scheme~(S) by the simpler
inference rule scheme
$X,Y\vd f$ for any objects $X$, $Y$
and any arrow $f\colon X\arr Y$ (due to p.~c) of Proposition~\ref{Proof1} above).

The situation with {\it coverings\/} of \Csub is not as simple: after all,
a covering in \Csub may contain simultaneously open arrows of {\it arbitrary high\/}
biheight, which indicates that there may exist no finite proof of such a covering in
the formal system with inference rules (C)--(PT4$''$),~(S).

Nevertheless, if one constructs first a {\it clopos\/}
$\la\Csub,\Osub\ra$
(where all proofs are finite), then the coverings
of \tsub have a simple description ``modulo
$\la\Csub,\Osub\ra$'' as shows the Corollary to the following proposition.
\begin{Prop}
\label{Proof2}
The inclusion functor $\Csub\subset\Shu$ reflects coverings.
In other words, a family $\{U_i\arr X\}_{i\in I}$ of open arrows
in \Csub is a covering in \Csub iff it is an epi family in \Shu.
\end{Prop}

\begin{Cor}
\label{cor3}
The inclusion functor $\Csub\subset\Shu$ reflects coverings.
In other words, a family $\{U_i\arr X\}_{i\in I}$ of open arrows
in \Csub is a covering in \Csub iff it is an epi family in \Shu.
\end{Cor}

\ToBeContinued
\fi
\section{Precloposes and Generalized Grothendieck Topologies}
\label{ggt}
A category \cC together with a class $\cO\subset\cC_1$ of arrows
of \cC will be called a {\it preclopos\/} iff \cO contains all
isomorphisms of \cC, is closed with respect to composition of arrows and,
besides, satisfies the following condition:

\noindent{\bf(qp)} For any commutative square
\begin{equation}
\label{sq}
\expandafter\ifx\csname xymatrix\endcsname\relax%
\divide\dgARROWLENGTH by2
\begin{diagram}
\node{Z}\arrow{e,t}{\alpha}\arrow{s,l}{\beta}\node{U}\arrow{s,r}{u}\\
\node{Y}\arrow{e,t}{f}\node{X}
\end{diagram}
\else
\xymatrix{
{Z}\ar[r]^{\alpha}\ar[d]_{\beta}&{U}\ar[d]^{u}\\
{Y}\ar[r]^{f}&{X}
}
\fi
\end{equation}
such that $u\in \cO$ there exists a commutative diagram
\begin{equation}
\label{sq1}
\expandafter\ifx\csname xymatrix\endcsname\relax%
\divide\dgARROWLENGTH by4
\dgDOTSPACING=0.15em
\begin{diagram}
\node{Z}\arrow[3]{e,t}{\alpha}\arrow[3]{s,l}{\beta}\arrow{se,..}
\node[3]{U}\arrow[3]{s,r}{u}\\
\node[2]{V}\arrow{ene,..}\arrow{ssw,r,..}{v}\\[2]
\node{Y}\arrow[3]{e,t}{f}\node[3]{X}
\end{diagram}
\else
\xymatrix@=8pt{
{Z}\ar[rrr]^{\alpha}\ar[ddd]_{\beta}\ar@{-->}[dr]&&&{U}\ar[ddd]^{u}\\
&V\ar@{-->}[urr]\ar@{-->}[ddl]_{v}\\
&&&\\
{Y}\ar[rrr]^{f}&&&{X}
}
\fi
\end{equation}
with $v\in\cO$.

Any clopos is a preclopos, evidently.

A sieve $R\subset X$ on an object $X$ of a preclopos \cC
will be said to be
\cO-{\it generated\/} or an \cO-{\it sieve\/} if there exists a family
$\{\sar U_i,u_i,X;\}_{i\in I}$ of arrows of \cO such that
the following condition is satisfied: an arrow $f$ belongs to $R$ iff
there exists $i\in I$ such that $f$ pulls through $u_i$:
\begin{equation}
\label{tri}
\expandafter\ifx\csname xymatrix\endcsname\relax%
\divide\dgARROWLENGTH by2
\dgDOTSPACING=0.15em
\begin{diagram}
\node[2]{U_i}\arrow{se,l,..}{u_i}\\
\node{Y}\arrow{ne,..}\arrow[2]{e,t,}{f}\node[2]{X}
\end{diagram}
\else
\xymatrix{
&{U_i}\ar@{-->}[dr]^{u_i}\\
{Y}\ar@{-->}[ur]\ar[rr]^{f}&&{X}
}
\fi
\end{equation}

\begin{Prop}
\label{pbsieve}
For any \cO-sieve $R\subset X$ and any arrow $f\colon Y\arr X$ the sieve
$f^*R$ is an \cO-sieve.
\end{Prop}

Proposition \ref{pbsieve} justifies the following definition.

A class $\tau=\bigcup\{\tau X:X\in\cC_0\}$ of \cO-sieves on a preclopos
$(\cC,\cO)$ will be called an \cO~-~{\it Gro\-then\-dieck topology\/}
on $(\cC,\cO)$ or, simply, an \cO-{\it topology\/} if:
\begin{list}{}{\topsep1mm\itemsep1mm}
\item[\bf(GT1)] For any object $X\in\cC_0$ the sieve $X\subset X$
belongs to $\tau$;
\item[\bf(GT2)] For any $R\in\tau X$ and any arrow $f\colon Y\arr X$
the sieve $f^*R$ belongs to $\tau Y$;
\item[\bf(GT3)] Let $R\in\tau X$ and $R'\subset X$ be an \cO-sieve.
If for any $v\colon Y\arr X$ in $R$ the sieve $v^*R'$ belongs to $\tau Y$,
then $R'$ belongs to $\tau X$ (local character).
\end{list}

\end{document}